\documentclass[3p,12pt]{elsarticle}
\usepackage{dsfont}
\usepackage{amssymb}
\usepackage{amsthm}
\usepackage{amsmath}
\usepackage{lineno}
\usepackage{empheq}
\usepackage{graphicx}
\usepackage{booktabs}
\usepackage{lineno}

\newdefinition{rmk}{Remark}
\newproof{pf}{Proof}
\newproof{pot}{Proof of Theorem \ref{thm2}}
\usepackage{caption}
\usepackage{subcaption}
\begin{document}
\begin{frontmatter}
\title{The stabilization of high-order multistep schemes for the Laguerre one-way wave equation solver.}
\author{Andrew V. Terekhov}
\ead{andrew.terekhov@mail.ru}
\address{Institute of Computational Mathematics and Mathematical Geophysics,
630090, Novosibirsk, Russia}
\begin{abstract}
This paper considers spectral-difference methods of a high-order of accuracy for solving the one-way wave equation using the Laguerre integral transform with respect to time as the base. In order to provide a high spatial accuracy and calculation stability, the Richardson method can be employed. However such an approach requires high computer costs, therefore we consider alternative algorithms based on the Adams multistep schemes. To reach the stability first for the 1D one-way equation and then for  the 2D case, the stabilizing procedures using the spline interpolation were developed. This made possible to efficiently implement a predictor-corrector type method in terms of which a boundary value problem for high-order elliptic equations is substituted for a sequence of inversions of second order elliptic operators thus decreasing computer costs. To assess the accuracy and stability of difference approximations for the 1D one-way wave equation, the analytical solution based on the double Laguerre transform was obtained. This solution can be efficiently calculated if for summation one makes use of fast algorithms of computing a discrete linear convolution. For the 2D one-way wave equation the stability and accuracy of the procedures proposed have been studied on implementing the migration algorithm within a problem of seismic prospecting.
\end{abstract}
\end{frontmatter}
\section{Introduction}
Mathematical models based on the one-way wave equation (OWWE) are often considered in problems of  ocean acoustics \cite{Tappert1977,Jensen2011,Lee2000}, seismic prospecting \cite{Claerbout1971,seis_interpr,Angus2013,Biondi2006} , as well as for setting non-reflecting boundary conditions \cite{Lindman1975,Engquist1977,Trefethen1986}
\begin{equation}
\frac{\partial \tilde{u}}{\partial z}=-\mathrm{i}\frac{\omega}{c}\sqrt{1-\left(\frac{c k_x}{\omega}\right)^2} \tilde{u},
\label{one-way}
\end{equation}
where $\mathrm{i}=\sqrt{-1}$, $\tilde{u}\equiv \tilde{u}(k_x,z,\omega)$ is a wave component at the angular frequency $\omega$, $k_x$ is the horizontal wave number, $c$ is the wave velocity, the vertical direction $z$ is the extrapolation direction, i.e., the direction of one-way propagation, and the positive axis $z$ is directed downward, i.e., toward increasing depth.
The square-root operator can be formally represented by the Pad\'e expansion \cite{Claerbout:1985,Halpern1988,Bamberger1988,Bamberger1988a}
\begin{equation}
\sqrt{1-\left(\frac{c k_x}{\omega}\right)^2} \approx\left[1-\sum_{s=1}^n\frac{\beta_sc^2 k_x^2}{\omega^2-\gamma_sc^2 k_x^2}\right],
\label{sqr-root-approx}
\end{equation}
where the coefficients $\gamma_s,\beta_s$ for the propagation angle should be optimized\cite{Halpern1988,Lee01101985}.
The velocity model was assumed to be homogeneous, although it also yields satisfactory results for inhomogeneous media. In the latter case this model correctly keeps kinematics of waves, but not their amplitudes. For a wide range of problems such an approximate one-way model is admissible as the correct account of amplitudes increases computer costs \cite{Zhang2005}.

The fundamental problem of the downward continuation algorithms of wave fields is the instability. If the coefficients  $\gamma_s,\beta_s$ are real, then for the angles around $\pi/2$ the argument of the square root becomes lesser than zero, the left-hand side of equation (\ref{sqr-root-approx}) being complex, while the right-hand side is still real, hence causing inconsistency in the approximation. This results in an improper propagation of the evanescent mode which should exponentially decay.  For stabilizing the real Pad\'e  approximation there are a few approaches \cite{Ristow01121994,Li1991,Bunks1995} that allow suppressing unstable components of a wave field. On the other hand, setting the coefficients  $\gamma_s,\beta_s$  to be complex \cite{Milinazzo1997,Lu1998,Yevick2000,Amazonas2007}, a better consistency of the right-hand and the left-hand sides of equation  (\ref{sqr-root-approx}) can be attained. From the physical viewpoint this means the introduction of artificial dissipation that restricts an increase in instability for evanescent waves. However the presence of strong gradients of the function of velocity of a medium, the use of the Marchuk-Strang type splitting for decreasing computer costs \cite{Marchuk1968,Marchuk1990,Strang1968} as well as the simulation of high-frequency wave fields, etc. can bring about both the instability of calculation and excess energy dissipation of waves to be simulated. This is explained by the fact that optimal values of the coefficients $\gamma_s,\beta_s$  are selected based on the frozen coefficients principle for a homogeneous medium, while calculations are carried out for inhomogeneous velocity models employing difference approximations and different decompositions for differential equations. Also, it is well to bear in mind that with the same number of terms in series  (\ref{sqr-root-approx})  the accuracy of the complex Pad\'e  approximation is somewhat lower than that real.

In \cite{Terekhov2017}, another approach was proposed  that provides the stability and high accuracy of calculations for the real coefficients $\gamma_s,\beta_s$. In this case the solution is sought for as a series in the Laguerre functions, while for increasing the spatial approximation accuracy the Richardson extrapolation\cite{Marchuk:99314} is used, which additionally restricts the instability of approximation (\ref{sqr-root-approx}). Earlier the stabilizing properties of the Richardson method for other problems were noted \cite{Gragg1965,Chinni1994}. However the Richardson method is not efficient because with its use it is needed to additionally solve the original problem on a mesh with a doubled number of nodes.  We will consider a stabilizing procedure based on the spline-filtration providing the stability of the Adams type multistep schemes of a high accuracy order \cite{Butcher2008}. This will allow carrying out calculations with the fifth order of accuracy instead of the fourth order thus to some extent decreasing computer costs.

In addition to the instability of the OWWE another fundamental problem is solving poorly conditioned systems of linear algebraic equations (SLAEs) in a frequency domain \cite{Ernst2012,Costa2013}.  As opposed to the Fourier transform, after employing the Laguerre transform the matrix of SLAEs is real and well-conditioned. Coefficients of the Laguerre series expansion have a recurrent dependence, therefore for their calculation it is required to solve the SLAEs with the same matrix and different right-hand sides, for example, with a parallel dichotomy algorithm \cite{terekhov:Dichotomy,Terekhov:2013,Terekhov2016}. Even a higher efficiency of calculations can be attained if one turns from solving difference problems for elliptic high-order operators to a sequence of problems for a second order operator. To this end, having stabilized the Adams schemes instability we will consider a method of the predictor-corrector type of the fifth accuracy order \cite{Butcher2008}.
\section{The stability Analysis for a model 1D one-way wave equation}
 The aspects of stability in constructing a numerical method for solving the 2D OWWE occupy a highly important place. To investigate the stability let us first consider a model problem for the 1D OWWE:

\begin{equation}
\label{advection}
\partial_t v+c\partial_x v=0,\quad t > 0,  \; x\in \mathbb{R}
\end{equation}
with the initial condition $v(x,0)=\varphi(x),\;(\varphi(0)=\varphi(1))$ and the periodic boundary condition $v(0,t)=v(1,t)$.

To solve problem  (\ref{advection}), let us consider the direct and inverse Laguerre transforms \cite{Integral_Transform} of a function $g(t)\in L_2(0,\infty)$
\begin{equation}
\label{series_lag}
\mathds{L}\{g(t)\}=\bar{g}_m=\int_{0}^{\infty}g(t)l_m(\eta t)dt,\quad\quad g(t)=\mathds{L}^{-1}\{\bar{g}_m\}=\sum_{m=0}^{\infty}\bar{g}_ml_m(\eta t),
\end{equation}
where $l_m(\eta t)\equiv\sqrt{\eta}\exp({-\eta t/2})L_m(\eta t)$ are the orthogonal Laguerre functions, $L_m(t)$ is the Laguerre polynomial of $m$ degree and $\eta>0$ is the transformation parameter.

Setting  $\mbox{$\lim_{t\rightarrow \infty}g(t)=0$}$, the following relations are valid \cite{Integral_Transform,Mikhailenko1999}
\begin{equation}
\mathds{L}\left\{\frac{d}{d t}g(t)\right\}=\frac{\eta}{2}\bar{g}_m+\Phi_1(\bar{g}_m),
\label{partial_t_lag}
\end{equation}
where
\begin{equation}
\Phi_1(\bar{g}_m)=\sqrt{\eta}g(0)+\eta\sum_{j=0}^{m-1}\bar{g}_j.
\label{phi_function}
\end{equation}

Making use of transformation (\ref{series_lag}) for equation (\ref{advection}), we obtain
\begin{equation}
\label{advection_laguerre}
\begin{array}{ll}
\displaystyle
\left(\frac{\eta}{2}+c\partial_x\right) \bar{v}^{m}+\Phi_1(\bar{v}^m)=0,\ m=0,1,...,
\end{array}
\end{equation}
where the index $m$ denotes number of a term in series (\ref{series_lag}).
Taking into consideration
$$\Phi_1(\bar{v}^{m})=\eta \bar{v}^{m-1}+\Phi_1(\bar{v}^{m-1}),$$
for studying the stability of difference schemes, let us turn to another form of equation  (\ref{advection_laguerre}):
\begin{subequations}
\label{advection_laguerre2}
  \begin{empheq}[left=\empheqlbrace]{align}
  \label{advection_laguerre2.a}
  \left(\frac{\eta}{2}+c\partial_x\right)\bar{v}^{0}+\sqrt{\eta}\varphi(x)=0,\\
  \label{advection_laguerre2.b}
  \left(\frac{\eta}{2}+c\partial_x\right)\bar{v}^{m}=\left(-\frac{\eta}{2}+c\partial_x\right)\bar{v}^{m-1},\ m=1,2,...
  \end{empheq}
\end{subequations}
For solving (\ref{advection_laguerre2.b}) let us consider the difference scheme of the first order of accuracy
\begin{equation}
\label{forward_scheme}
\frac{c\left(\bar{v}^{m}_{j+1}-\bar{v}^{m}_{j}\right)}{h_x}+\frac{\eta}{2}\bar{v}^{m}_{j}=\frac{c\left(\bar{v}^{m-1}_{j+1}
-\bar{v}^{m-1}_{j}\right)}{h_x}-\frac{\eta}{2}\bar{v}^{m-1}_{j}.
\end{equation}

Substituting the solution in the form
$
\bar{v}^{m}_{j}=\tilde{\bar{v}}^{m}\exp(\mathrm{i}k_xjh_x)
$
into difference equation (\ref{forward_scheme}), obtain
\begin{equation}
\tilde{\bar{v}}^{m}=\frac{\exp(\mathrm{i}k_xh_x)-\beta-1}{\exp(\mathrm{i}k_xh_x)+\beta-1}\tilde{\bar{v}}^{m-1}=G(k_x)\tilde{\bar{v}}^{m-1}.
\label{stability_forward}
\end{equation}
Here $\beta=\eta h_x/(2c)$, $G=G(k_x)$ is called the amplification factor which is a complex function of the wavenumber $k_x$. A difference equation will be stable in the Von Neumann sense \cite{CHARNEY1950} if $|G(k_x)|\leq1$  $\forall$ $k_x$. For equation (\ref{stability_forward}) let us estimate the value
$$
\left|G(k_x)\right|^2={\frac { \left( {\beta}+1-\cos \left( {k_xh_x} \right)  \right)^{2}+\sin ^{2} \left( {k_xh_x} \right) }{ \left( {
\beta}-1+\cos \left( {k_xh_x} \right)  \right) ^{2}+ \sin ^{2}
\left( {k_xh_x} \right)  }}=\frac{A}{B}.
$$
For  $c>0$ obtain $A>0,B>0,A-B=4\beta\left(1-\cos(k_xh_x)\right)\geq0$, hence, $\left|G(k_x)\right|^2\geq1$ and scheme (\ref{forward_scheme}) will be unstable. For  $c<0$ it can be shown that $\left|G(k_x)\right|^2\leq1$, which suggests the stability of the scheme.

Now let us consider another method of the first order of accuracy
\begin{equation}
\label{explicit_euler}
\frac{c\left(\bar{v}^{m}_{j+1}-\bar{v}^{m}_{j}\right)}{h_x}+\frac{\eta}{2}\bar{v}^{m}_{j+1}=
\frac{c\left(\bar{v}^{m-1}_{j+1}-\bar{v}^{m-1}_{j}\right)}{h_x}-\frac{\eta}{2}\bar{v}^{m-1}_{j+1}.
\end{equation}
In a similar manner reducing scheme (\ref{explicit_euler})  to the form  \mbox{$\tilde{\bar{v}}^{m}=G(k_x)\tilde{\bar{v}}^{m-1}$}, we obtain that for $c>0$ the value $A>0,B>0,A-B=4\,{\beta}\left(\cos \left( {k_xh_x} \right)-1  \right)\leq0$, and hence, $\left|G(k_x)\right|^2\leq1$. Thus, scheme (\ref{explicit_euler}) will be stable for $c>0$ and unstable for $c<0$.

Let us consider the Crank-Nicolson scheme (CN-Scheme) \cite{CrankNicolson2} of the second order of accuracy
\begin{equation}
\label{crank-nicolson}
\frac{c\left(\bar{v}^{m}_{j+1}-\bar{v}^{m}_{j}\right)}{h_x}+\frac{\eta}{4}\left(\bar{v}^{m}_{j+1}+\bar{v}^{m}_{j}\right)
=\frac{c\left(\bar{v}^{m-1}_{j+1}-\bar{v}^{m-1}_{j}\right)}{h_x}-\frac{\eta}{4}\left(\bar{v}^{m-1}_{j+1}+\bar{v}^{m-1}_{j}\right).
\end{equation}
It is not difficult to obtain that for the scheme in question $|G(k_x)|=1$ holds, therefore scheme (\ref{crank-nicolson}) is unconditionally stable.

For multistep schemes of a high-order of accuracy of the Adams-Moulton type (AM scheme) the stability has been analytically studied. Having omitted cumbersome operations we will restrict ourselves to the main conclusions. The AM-schemes of the third and fourth orders of accuracy will be unstable for periodic boundary conditions if $c>0$ and stable if $c<0$. The schemes of the fifth order of accuracy and higher are absolutely unstable. The Von Neumann spectral feature is a necessary but not sufficient condition of stability \cite{Godunov_Rib}, and in the case of non-periodic boundary conditions the stability is not preserved for the AM-schemes of higher than the second order of approximation.
To overcome these difficulties we will consider the ways of stabilizing the multistep AM-schemes of a high-order of accuracy first for the 1D and then for the 2D OWWE.
\section{The stabilization of  high-order schemes for the 1D one-way wave equation}
For equation (\ref{advection}) at  $c>0$  instead of the periodic boundary conditions we consider initial and boundary conditions of the form:
\begin{equation}
\label{boundary2}
\begin{array}{ll}
v(0,t)=f(t), & t\geq 0,\\
v(x,0)=0, & x\geq 0,\\
f(0)=0.
\end{array}
\end{equation}
A stable spectral-difference algorithm for solving the 2D OWWE was proposed in \cite{Terekhov2017}. It includes the Richardson extrapolation procedure that for problem (\ref{advection_laguerre}), (\ref{boundary2})  and scheme (\ref{crank-nicolson}) can be written down in the following form.

\textbf{Algorithm 1-1. The Richardson extrapolation.}

Let the auxiliary functions $\bar{v}^m(\Omega_1)$, $\bar{v}^m(\Omega_2)$ be defined on the meshes $\Omega_1,\Omega_{2}$ with the mesh steps $h_x$ and $h_x/2$.
To calculate the functions $\bar{v}^m$ accurate to $O(h_x^4)$, the following is necessary:
 \begin{enumerate}
   \item Based on the cubic splines interpolate values of the function $\Phi_1(\bar{v}^m)$, preset on the mesh $\Omega_1$, into nodes $\Omega_2$.
   \item On the mesh $\Omega_1$, applying equation (\ref{crank-nicolson}), calculate the solution $\bar{v}^{m}(\Omega_1)$.
   \item On the mesh $\Omega_2$, applying equation (\ref{crank-nicolson}), calculate the solution $\bar{v}^{m}(\Omega_2)$.
   \item Based on the Richardson extrapolation, correct the mesh function with $$\bar{v}^m=\frac{1}{3}\left(4\bar{v}^m(\Omega_2)-\bar{v}^m(\Omega_1)\right).$$
    \item Turn to the calculation of the  $(m+1)\text{th},\;\text{the }(m+2)\text{th}$, etc. coefficients of the expansion of the Laguerre series.
 \end{enumerate}
This technique of calculating the Laguerre series coefficients is stable and provides the fourth order of accuracy. However the necessity of calculating $\bar{v}^m(\Omega_2)$  triples the common computer costs, therefore there arises a problem of constructing a more efficient method of no less than fourth order of accuracy.

Let us consider a difference approximation for equation (\ref{advection_laguerre}) based on the Adams-Moulton multistep method of the fifth order of accuracy
\begin{equation}
\label{adams5}
\begin{array}{ll}\displaystyle
c\frac{\bar{v}^m_{i+1}-\bar{v}^m_{i}}{h_x}
=-\frac{1}{720}\sum_{j=-3}^1\alpha_j\left(\frac{\eta}{2}\bar{v}^m_{i+j}+\Phi_1\left(\bar{v}^m_{i+j}\right)\right),
\end{array}
\end{equation}
where the coefficients of the difference scheme are equal to $\alpha_{-3}=-19,\;\alpha_{-2}=106,\;\alpha_{-1}=-264,\;\alpha_0=646,\;\alpha_1=251$.
This scheme is unstable according to the Von Neumann spectral property, but it can be stabilized when carrying out calculations in the following manner.

\textbf{Algorithm 1-2. The stabilization of the Adams-Moulton scheme via the quintic spline filtration.}

To calculate the functions $\bar{v}^m$ accurate to $O(h_x^5)$ on the mesh $\Omega$  with the mesh step $h_x$, the following is necessary:
 \begin{enumerate}
\item  Let the number of nodes of the mesh $\Omega$ be odd. Construct the quintic splines for the function $\Phi_1({\bar{v}}^m)$ using only odd nodes of the mesh.
\item  Replace values of the function  $\Phi_1(\bar{v}^m_{k})$ for even  $k$ by their interpolated values (the quintic spline filtration).
\item  Applying equation (\ref{adams5}), calculate the solution $\bar{v}^{m}$.
\item  Turn to the calculation of the $(m+1)\text{th},\;\text{the} (m+2)\text{th}$, etc. coefficients of the expansion of the Laguerre series.
 \end{enumerate}
Such an algorithm of calculations makes possible to stabilize the numerical instability of scheme  (\ref{adams5}) and to attain a higher approximation order as compared to the Richardson extrapolation. To stabilize the solution, instead of the quintic spline one can use other interpolation algorithms \cite{Marchuk1982,Berrut2004}:
the cubic spline interpolation, barycentric, Lagrangian, etc. However numerous computer-aided experiments have not revealed any advantages over splines because the procedures of constructing splines are efficient enough as compared to solving elliptic equations in the 2D case. In addition, the barycentric interpolation demands high computer costs and both the Lagrangian interpolation and the cubic spline-interpolation are more dissipative than the quintic splines. If splines are not being used, the nodes of an interpolating polynomial should be symmetrically placed regarding the node for which the interpolated value is calculated. Otherwise due to the asymmetry of interpolation nodes the profile of a wave is distorted or the instability of calculation arises.

We can offer another way of stabilizing the numerical instability of the  AM-schemes, which does not demand the calculation of splines.

\textbf{Algorithm 1-3. The stabilization of the Adams-Moulton scheme via inconsistent approximation.}
\begin{enumerate}
\item
Taking into account the equivalence of problems (\ref{advection_laguerre}) and (\ref{advection_laguerre2}), for computing the values of the grid functions $\Phi_1\left(\bar{v}^m_{i}\right)$ instead of  (\ref{phi_function}) use the following approximation
\begin{equation}
\Phi_1\left(\bar{v}^m_{i}\right)=-\frac{\eta}{2}\bar{v}^{m-1}_i+c\frac{-\bar{v}^{m-1}_{i+2}+
8\bar{v}^{m-1}_{i+1}-8\bar{v}^{m-1}_{i-1}+\bar{v}^{m-1}_{i-2}}{12h_x}+O\left(h_x^4\right).
\label{AM5-D4}
\end{equation}
\item Solve equation (\ref{advection_laguerre}) through scheme (\ref{adams5}).
\item Turn to the calculation of the $(m+1)\text{th},\;\text{the }(m+2)\text{th}$, etc. coefficients of the expansion of the Laguerre series.
\end{enumerate}
The non-consistent approximation of the operator $\partial/\partial x$ for the right-and the left-hand sides of equations (\ref{advection_laguerre}),(\ref{advection_laguerre2}), stipulates supplementary non-physical dissipation preventing the development of instability. However, if instead of the central approximation of the fourth order of accuracy in (\ref{AM5-D4}) one uses a higher order approximation or a non-central scheme, the stability is lost.

In addition to the above we have considered a stable algorithm based on the central differences of a high-order of accuracy for the operator  $\partial/\partial x$. For the 1D OWWE such calculations were carried out up to the schemes accurate to the twelfth order. As a result, the stability and high-order of accuracy were attained. However this approach will not work for the 2D OWWE because of high computer costs that are required for solving SLAEs.

Thus, in addition to the Richardson method, other stable algorithms of a high-order of accuracy can be proposed. However for calculating of  $\bar{v}^{m}$   the function $\bar{v}^{m-1}$ should be known throughout the whole calculation domain as its values are needed for implementing the stabilizing procedures.

\section{The analytic solution via the Laguerre transforms for the 1D one-way wave equation}
In order to assess the accuracy of the algorithms proposed, let us consider a fully analytical method for solving the 1D OWWE. To satisfy boundary conditions (\ref{boundary2}), we seek the solution to equation (\ref{advection_laguerre2})  in the form
\begin{equation}
\bar{v}^m(x)=\sum_{j=0}^{\infty}V^m_jl_j(\kappa x),\quad m=0,1,2...,
\label{laguerre_x_expansion}
\end{equation}
where the transformation parameter $\kappa>0$. Then, after applying the Laguerre spatial transform to equation  (\ref{advection_laguerre2}) with allowance for initial boundary conditions  (\ref{boundary2}) we have
\begin{subequations}
\label{lagx}
\begin{empheq}[left=\empheqlbrace]{align}
\label{lagx.1}
\displaystyle \left(\eta+c\kappa\right)V^m_0=\left(-\eta+c\kappa\right)V^{m-1}_0+2c\sqrt{\kappa}\left(\bar{f}^m- \bar{f}^{m-1}\right),&\quad m=0,1,...,\\
\label{lagx.2}
\left(\eta+c\kappa\right)V^m_j+2c\Upsilon(V^m_j)=\left(-\eta+c\kappa\right)V^{m-1}_{j}+2c\Upsilon(V^{m-1}_{j}),&\quad m=0,1,...;\, j=1,2,..,
\end{empheq}
\end{subequations}
where
\begin{equation}
\Upsilon\left(V^m_j\right)=\kappa \sum_{i=0}^{j-1}{V^m_i}=\kappa V^m_{j-1}+\Upsilon\left(V^m_{j-1}\right),
\label{eq1}
\end{equation}
$$V^{m}_j\equiv0,\;\bar{f}^{m}\equiv0,\quad \forall \; m<0.$$
Taking (\ref{eq1}) into account, equation (\ref{lagx.2}) takes the following form
\begin{equation}
\left(\eta+c\kappa\right)V^m_j+\left(\eta-c\kappa\right){V}^{m-1}_j=
\left(\eta-c\kappa\right)V^{m}_{j-1}+\left(\eta+c\kappa\right)V^{m-1}_{j-1},\quad  m=0,1,...;\,j=1,2,...
\end{equation}
Since  $c>0$, then selecting $\kappa=\eta/c$, we finally obtain
\begin{equation}
\label{lagx3}
\left\{\begin{array}{ll}
V^m_0=\kappa^{-1/2}\left(\bar{f}^m-\bar{f}^{m-1}\right),& m=0,1,...,
\\
V^m_j=V^{m-1}_{j-1},& m=0,1,...;\; j=1,2,...
\end{array}
\right.
\end{equation}
Based on (\ref{laguerre_x_expansion}),(\ref{lagx3}), we can write down the solution to equation   (\ref{advection_laguerre2}),(\ref{boundary2}) in the form
\begin{equation}
\bar{v}^m(x)=\sum_{j=0}^{\infty}V^m_jl_j(\kappa x)=\sum_{j=0}^mV^{m-j}_{0}l_{j}(\kappa x),\quad m=0,1,2,...
\label{conv}
\end{equation}
The latter sum in (\ref{conv}) is a discrete linear convolution, therefore for a given $x$ the functions $\bar{v}^m(x),\; m=0,...,M$ can be computed in $O(M\log M)$ arithmetical operations based on the algorithm of the FFT \cite{Nussbaumer1982}. Note that if we select $\kappa\neq \eta/c$, the solution for $\bar{v}^m$ will not be representative in the form of convolution thus increasing computer costs. The final solution to equation (\ref{advection}) in the time domain is calculated via the inverse Laguerre transform (\ref{series_lag}).
\section{Numerical Experiments for the 1D one-way wave equation}
\label{sec:1dowwe}
For testing  the methods proposed for the 1D OWWE, we used a homogeneous medium model with the speed $3000\;m/s$ and the size $7.5\;km$. For the test calculations we set  boundary condition (\ref{boundary2}) depending on time as
\begin{equation}
f(t)=\exp\left[-\frac{(2\pi f_0(t-t_0))^2}{\delta^2}\right]\sin(2\pi f_0(t-t_0)),
\label{source}
\end{equation}
where $t_0=0.2s,\,\delta=4,\,f_0=30\text{Hz}$. As compared to the Fourier transform, where the basis functions are uniquely defined, the parameter $\eta$  for using the Laguerre transform (\ref{series_lag}) should be set. This parameter was experimentally chosen based on the analysis of the convergence rate of the Fourier-Laguerre series for the shifted function $f(t)$ with $t_0=T$, where T is the upper boundary of the time interval for which the wave field is calculated. The parameter $\eta$ is chosen such that the function f(t) with $t_0 = T$ in the mean-quadratic norm is approximated accurate to $\varepsilon < 10^{-10}$. The number of addends in series (\ref{series_lag}) was $n=2500$ for $T=2\;s$; the expansion parameters were $\eta=600$.

From Table~\ref{tabl:advection} it is evident that with decreasing the mesh size by a factor of two, the error of the Adams and the Richardson methods is decreasing according to theoretical approximation order. For example, for meshes with the number of nodes $N_x=2000$ and  $N_x=4000$ the values of the error of the AM-scheme of fifth order with the quintic spline interpolation (AM5-I5) has $31$ times  difference, which almost corresponds to the fifth order of approximation. The AM-scheme of the sixth order with the seventh order spline interpolation (AM6-I7) demonstrates the sixth order of approximation, while the AM-scheme of the fifth order with formula (\ref{AM5-D4}) (AM5-D4) is only of the fourth order of approximation.

\begin{figure}[!htb]
\centering

        \begin{subfigure}[b]{0.49\textwidth}
                \includegraphics[width=\textwidth]{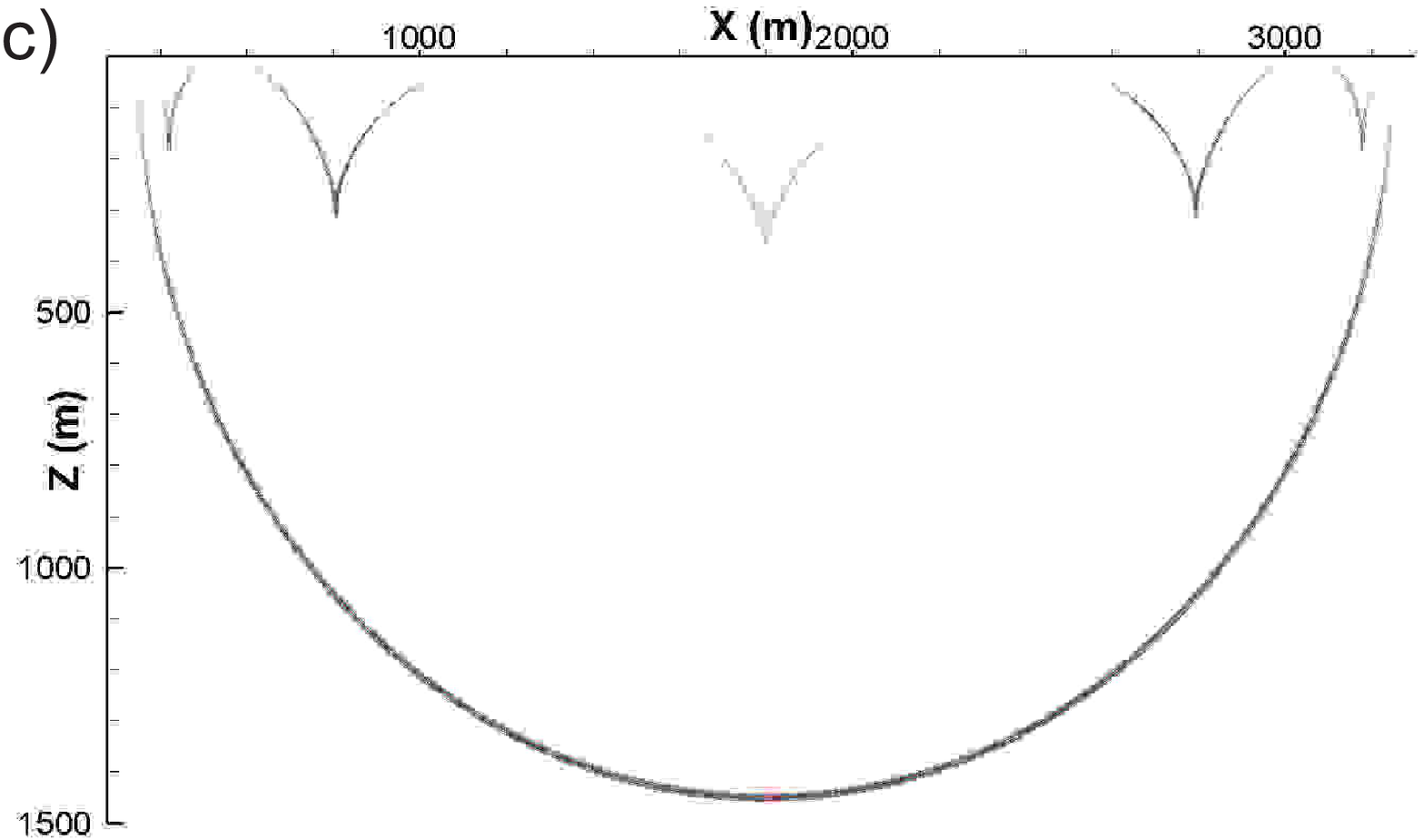}
                \label{fig:homo21}
        \end{subfigure}
             ~
                \begin{subfigure}[b]{0.49\textwidth}
                \includegraphics[width=\textwidth]{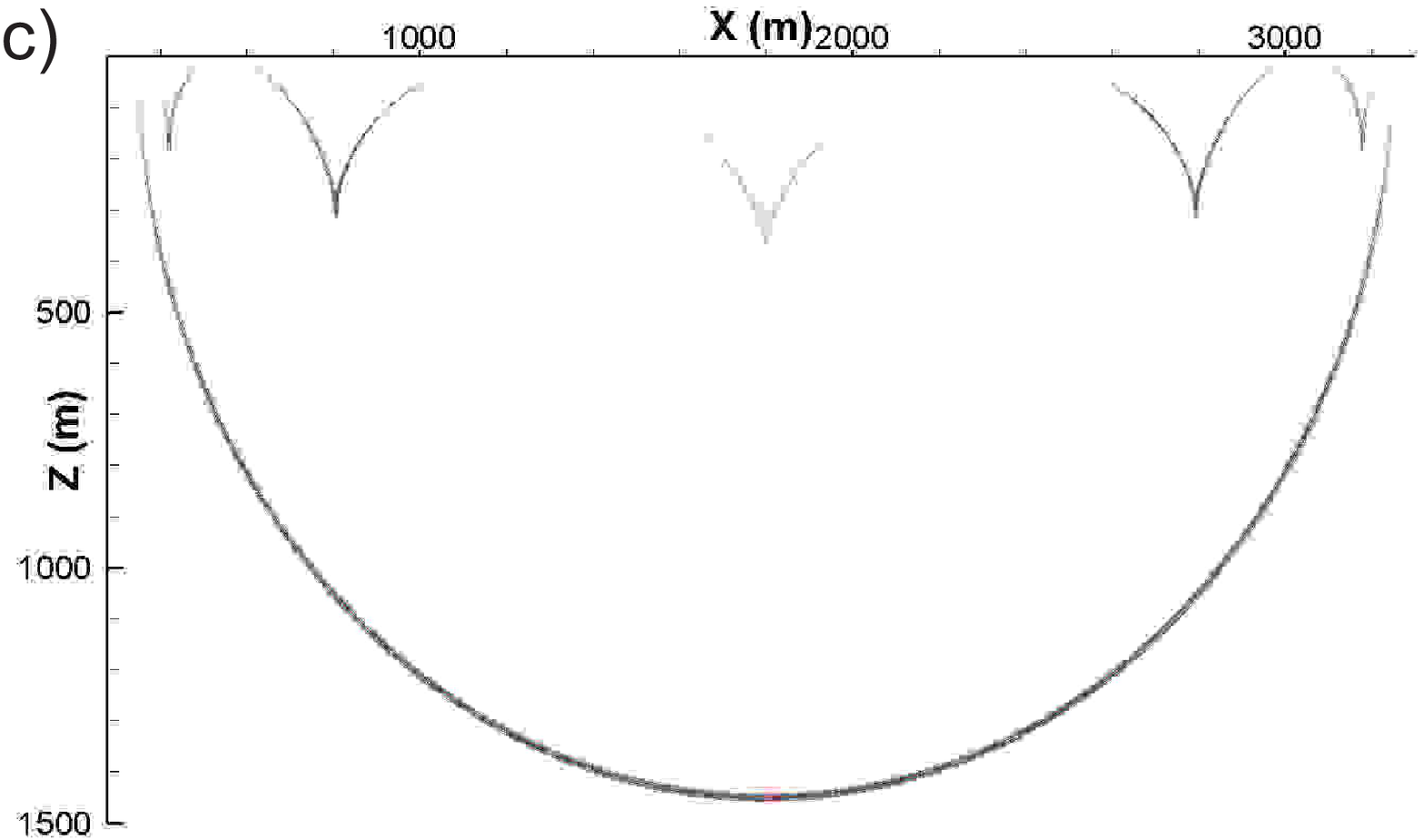}
                \label{fig:homo22}
        \end{subfigure}%

\caption{Dependence of the wave field on the coordinate for the method proposed and  different meshes (a) $N_x=1000$, (b) $N_x=2000$. }
\label{fig:advection1}
\end{figure}

From Table~\ref{tabl:advection} and Fig.~\ref{fig:advection1} it also follows that with an equal mesh size, the Richardson extrapolation of the fourth order of approximation is more accurate as compared to the Adams methods of the fifth and sixth orders. There is no contradiction because the estimation of the accuracy of difference schemes  includes a constant independent of the mesh size. This constant is smaller for the Richardson method as compared to the AM5-15 method as all splines are constructed on a mesh with a doubled step. In addition, the Richardson method requires the solution to a supplementary problem on the mesh $\Omega_2$, therefore it is more correct to compare the accuracy of calculations when the general number of nodes of the meshes $\Omega_1,\Omega_2$ is equal  to the number of nodes of the mesh $\Omega$ for other methods. Indeed, comparing solutions for the Richardson method with  $N_x=1500$ for the mesh $\Omega_1$ and the Adams methods with $N_x=4500$, it is evident  (Table ~\ref{tabl:advection}) that the latter are significantly more accurate. Applying the spline-filtration procedure does not bring about a considerable loss in accuracy, otherwise the accuracy of the methods  AM5-I5, AM6-I7 would be lower or compatible with the fourth order method AM5-D4,  for which the filtration is not employed.

For comparison similar calculations were carried out for the CN-Scheme and the explicit Runge-Kutta method (the RK4-method) \cite{Butcher2008}, which are of the second and fourth orders of accuracy, respectively. For the method RK4, values of the grid function $\Phi_1\left(\bar{v}^m\right)$ in semi-integer nodes of the mesh were calculated based on the quintic splines. The implementation of the implicit Runge-Kutta method of a high-order is not reasonable because of essentially higher computer costs as compared to the approaches proposed in the given study. From Table~\ref{tabl:advection} and Figure~\ref{fig:advection1} it is clear that the method RK4 and the CN-scheme are stable  and converge to the analytical solution with decreasing a mesh size.  However for large mesh steps, the method RK4 possesses a pronounced numerical dissipation, while on the contrary, the Crank-Nicolson scheme demonstrates a dispersive error. A low accuracy for large mesh steps makes the application of these algorithms disadvantageous in comparison with the Adams methods.
\begin{table}
\center \small
\begin{tabular}{lcccccc}
  \hline
   & AM5-I5&AM6-I7&CN&RK4 & Richardson&AM5-D4 \\
  $N_x$&&&&&&   \\\hline
 1000& 0.32& 0.18  & 1.47 & 0.99 & 6.04e-02&0.56 \\
 1500&6.67e-2 &2.4e-2& 1.51 & 0.92 & 1.13e-2 & 0.18 \\
 2000&1.72e-2  &4.6e-3&1.38 & 0.6 & 3.5e-3 & 6.5e-2 \\
 3000&2.3e-3 &4.18e-4& 0.87 & 1.68 & 6.82e-4 & 1.33e-2 \\
 4000&5.6e-4 &7.52e-5& 0.53 & 5.46e-2 & 2.14e-4 &4.2e-3 \\
 4500&3.1e-4 &3.72e-5& 0.43 & 3.39e-2 & 1.33e-4 &2.6e-3 \\\hline
\end{tabular}
\caption{Dependence of the error value  $\|u^{\mathrm{exact}}-u^h\|_{2}/\|u^{\mathrm{exact}}\|_{2}$ on the number of mesh nodes for different methods.}
\label{tabl:advection}
\end{table}

To evaluate dissipative properties of the methods proposed, let us consider the integral of the form
\begin{equation}
K(x)=\int_{0}^{\infty} v^2(x,t) dt=\sum_{k=0}^{\infty}\left[\bar{v}^{k}(x)\right]^2,
\label{pseudo_energy}
\end{equation}
where the latter equality is the Parseval relation.

\begin{figure}[!htb]
\centering
        \begin{subfigure}[b]{0.49\textwidth}
                \includegraphics[width=\textwidth]{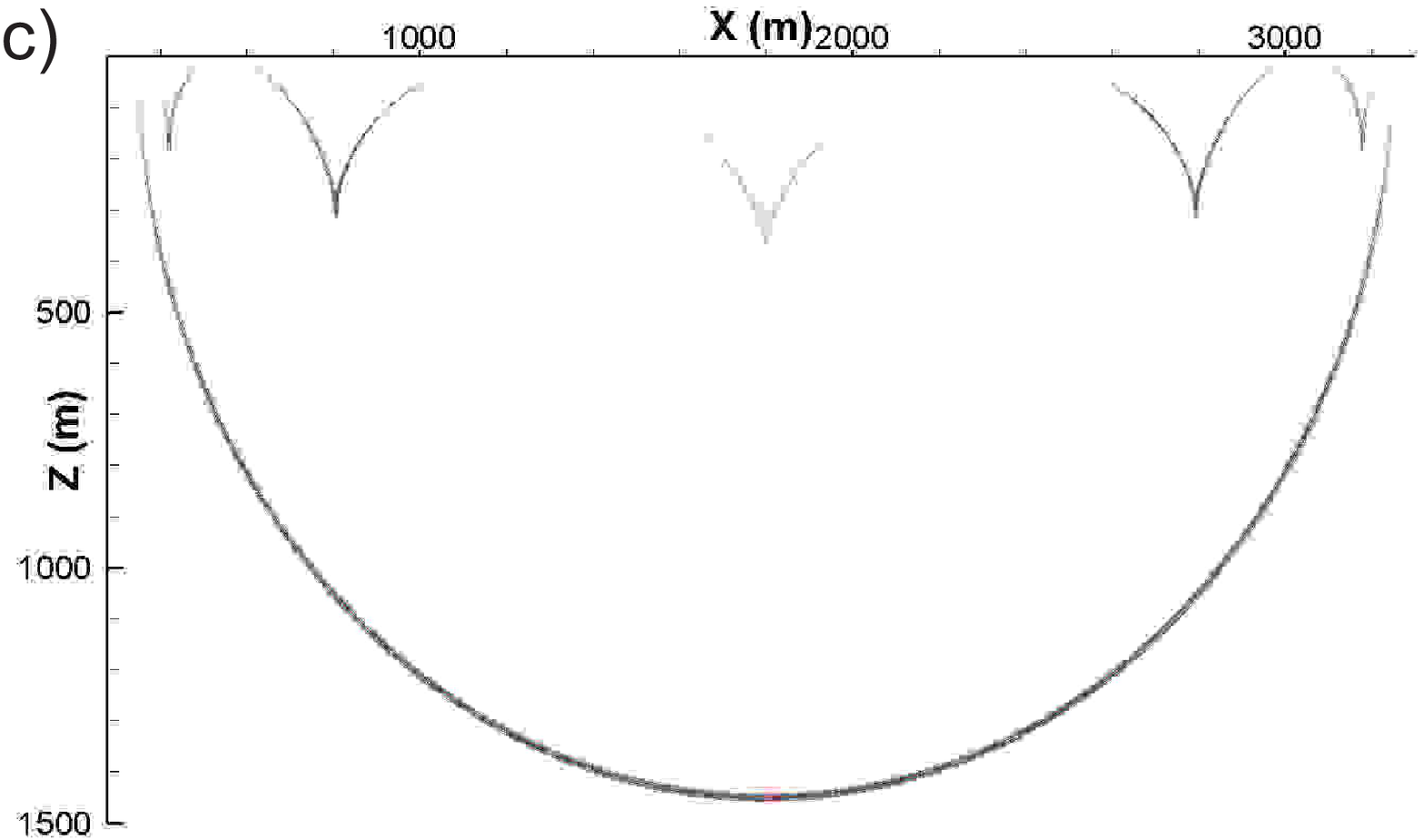}
                \label{fig:homo21}
        \end{subfigure}
        ~
        \begin{subfigure}[b]{0.49\textwidth}
                \includegraphics[width=\textwidth]{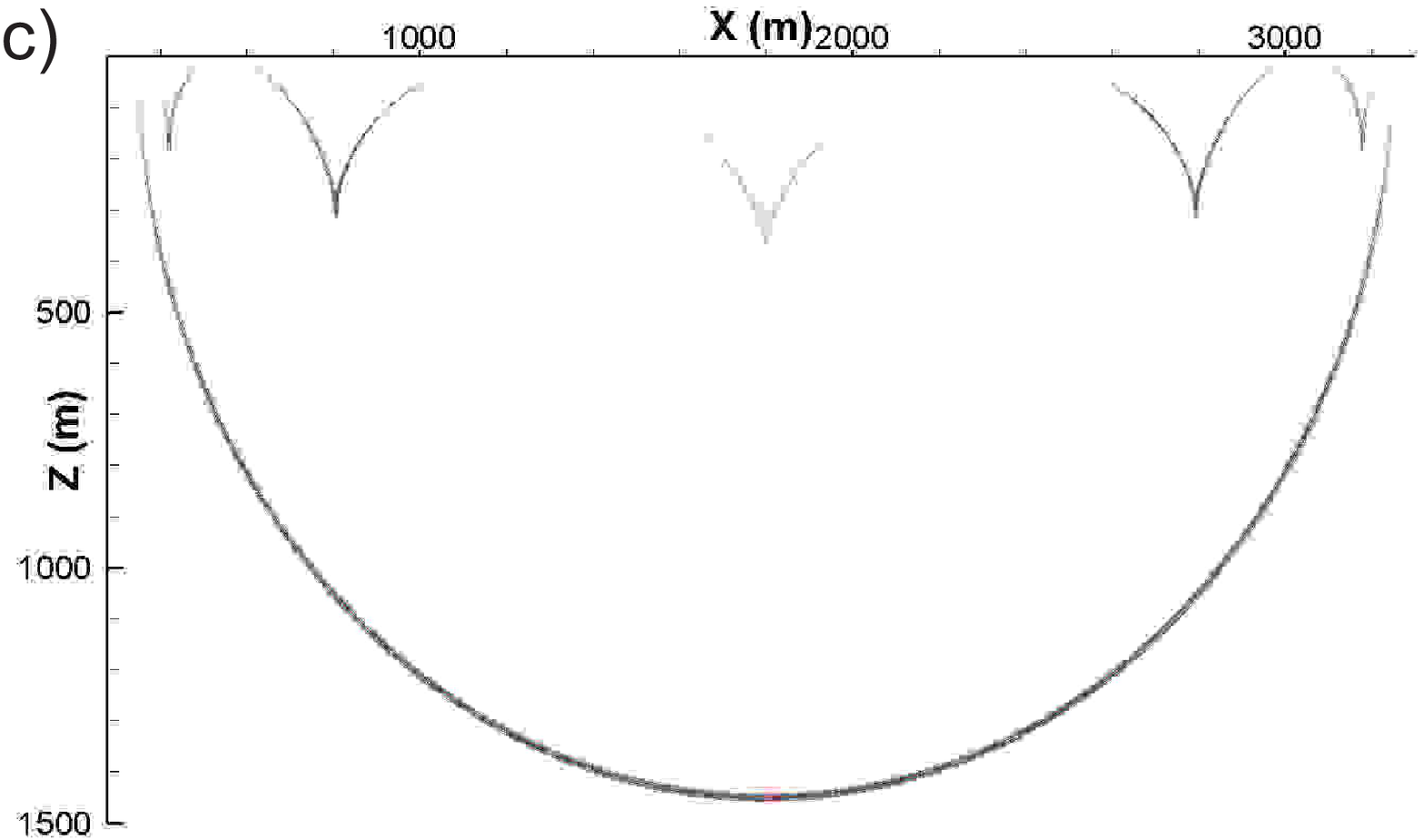}
                \label{fig:homo21}
        \end{subfigure}

\caption{Dependence of the value $K(x)$ on the coordinate for   different methods and meshes: (a) $N_x=1000$, (b) $N_x=2000$.}
                \label{fig_u2}
\end{figure}
For problem (\ref{advection}),(\ref{boundary2})  at $c=const>0$  with a sufficient number of terms in  Laguerre series (\ref{series_lag}), $K(x)=const$ should hold with a good accuracy. Figure~\ref{fig_u2}a,b shows that for the analytical method  (\ref{lagx3}),(\ref{conv})  and for the CN-scheme the value $K(x)$ is preserved with the precision of a machine. However due to the numerical dispersion the solution obtained with the CN-scheme does not satisfy the 1D OWWE equation with some kind of accuracy for large mesh steps. The Richardson method at $N_x=1000$  (Fig.~\ref{fig_u2}a) is less dissipative than the algorithms AM5-I5 and AM6-I7, while with an increase of the number of mesh nodes  (Fig.~\ref{fig_u2}b) the situation is contrary. This means that when solving the 2D OWWE by the Richardson method the stability will be stronger as compared to the Adams methods.
The method RK4 is most dissipative among all under consideration, hence the initial impulse from a source has smoothed into the straight line (Fig.~\ref{fig:advection1}a). Thus, the explicit method  RK4 and the CN-scheme, as was already mentioned, cannot be offered for the use within the Laguerre method.

\section{The 2D one-way wave equation solver}
\subsection{Temporal approximation}
Let us write down 2D OWWE (\ref{one-way}),(\ref{sqr-root-approx}) for the spatial-temporal domain \cite{Bamberger1988,Bamberger1988a}
\begin{subequations}
\label{main_eq}
\begin{empheq}[left=\empheqlbrace]{align}
\label{main_eq1}
\frac{\partial u}{\partial t}+c\frac{\partial u}{\partial z}-\sum_{s=1}^{n}\frac{\partial \psi_s}{\partial t}=0,\\
\frac{1}{c^2}\frac{\partial^2 \psi_s}{\partial t^2}-\gamma_s \frac{\partial^2\psi_s }{\partial x^2}-\beta_s\frac{\partial^2 u }{\partial x^2}=0,\quad s=1,2,...n,
\label{main_eq2}
\end{empheq}
\end{subequations}
where $u\equiv u(x,z,t)$ is the field variable, $\psi_s\equiv \psi_s(x,z,t)$ are auxiliary functions.

Assuming $g(0)=\left.\frac{dg}{dt}(t)\right|_{t=0}=0$ and $\lim_{t\rightarrow \infty}g(t)=\lim_{t\rightarrow \infty}\frac{dg}{dt}(t)=0$, we can show \cite{Integral_Transform,Mikhailenko1999} that
\begin{equation*}
L\left\{\frac{d^{2}}{d t^2}g(t)\right\}= \left(\frac{\eta}{2}\right)^{2}\bar{g}_m+\Phi_2(\bar{g}_m),\quad\Phi_2(\bar{g}_m)\equiv\eta^2\sum_{j=0}^{m-1}(m-j)\bar{g}_j.
\label{partial_t_lag}
\end{equation*}
Then applying the Laguerre transform (\ref{series_lag})  to equations (\ref{main_eq})  we obtain the following system of equations for the calculation of the $m$-th coefficient of expansion:
\begin{subequations}
\label{main_eq_laguerre}
\begin{empheq}[left=\empheqlbrace]{align}
\label{main_eq1_laguerre}
\tilde{\eta} \bar{u}^m+c\frac{\partial \bar{u}^m}{\partial z}=\sum_{s=1}^3\left(\tilde{\eta}\bar{\psi}_s^m+\Phi_1\left(\bar{\psi}_s^m\right)\right)-\Phi_1\left(\bar{u}^m\right),\\
\label{main_eq2_laguerre}
c^2\gamma_s\frac{\partial^2\bar{\psi}_s^m }{\partial x^2}-\tilde{\eta}^2\bar{\psi}_s^m+\beta_sc^2\frac{\partial^2\bar{u}^m }{\partial x^2}=\Phi_2(\bar{\psi}_s^m) ,\quad s=1,2,3,
\end{empheq}
\end{subequations}
where $\tilde{\eta}=\eta/2$ and the index $m$ denotes number of a term in  series (\ref{series_lag}). The polynomial coefficients $\gamma_s,\beta_s$ for $n=3$ are chosen as follows: ~$\gamma_1=0.972926132$, $\gamma_2=0.744418059$, $\gamma_3=0.150843924$, $\beta_1=0.004210420$, $\beta_2=0.081312882$, $\beta_3=0.414236605$, for which, as shown in \cite{seis_interpr,Lee01101985}, such approximation is valid up to the angles of $89$ degrees.
\subsection{The spatial approximation}
Multistep schemes of the Adams type of a high-order for solving  the 2D OWWE are of practical importance as opposed to the Richardson method which requires solving an auxiliary problem on the mesh  $\Omega_2$   with a doubled number of nodes.
To approximate equation (\ref{main_eq_laguerre}) we will use the Adams-Moulton scheme of the fifth order of accuracy:
\begin{subequations}
\label{main_eq_diff}
\begin{empheq}[left=\empheqlbrace]{align}
\frac{\bar{u}^m_{ik+1}-\bar{u}^m_{ik}}{h_z}=\frac{1}{720c}\sum_{j=-3}^1\alpha_j\left(\sum_{s=1}^{3}\left(\tilde{\eta}\bar{\psi}^{m,s}_{ik+j}
+\Phi_1\left(\bar{\psi}^{m,s}_{ik+j}\right)\right)-\tilde{\eta}\bar{u}^m_{ik+j}-\Phi_1\left(\bar{u}^m_{ik+j}\right)\right)
\label{adams_owwe.a}
,\\
\label{adams_owwe.b}
c^2\gamma_s \mathcal{L}_{x}\bar{\psi}^{m,s}_{ik+1}-\tilde{ \eta}^2\bar{\psi}^{m,s}_{ik+1}=-c^2\beta_s \mathcal{L}_{x}\bar{u}^m_{ik+1}+\Phi_2(\bar{\psi}^{m,s}_{ik+1}),\quad s=1,2,3,
\end{empheq}
\end{subequations}
where the coefficients of a difference scheme $\alpha_{-3}=-19,\;\alpha_{-2}=106,\;\alpha_{-1}=-264,\;\alpha_0=646,\;\alpha_1=251$ and
the difference operator $\mathcal{L}_{x}$ is of the form
\begin{equation}
\mathcal{L}_{x}f(x)\equiv\frac{1}{h_x^2}\left[a_0f(x)+\sum_{j=1}^{N}a_j\left(f(x-jh_x)+f(x+jh_x)\right)\right]=\frac{\partial^2 f}{\partial x^2}(x)+O(h_x^{ 2N}).
\label{drp_coeff}
\end{equation}
For approximating $\partial^2/\partial x^2$  it is reasonable to use the dispersion-relationship-preserving method (DRP) by Tam and Webb \cite{Tam1993262}, for which according to the Fourier derivative rule, $k_j\Longleftrightarrow-\mathrm{i}\partial_j,$
values of the optimized coefficients $a_n$ in (\ref{drp_coeff}) are defined as solution to the problem of minimizing the error functional in the space of wave numbers.
This approach and its various modifications \cite{GPR:GPR12160,Chu2012,Zhang2013511} make possible to decrease the number of mesh nodes and to preserve high accuracy of calculations as compared to conventional difference schemes obtained with the Taylor expansion in series \cite{Marchuk1982}. To provide the twelfth approximation order the coefficients of difference scheme (\ref{drp_coeff}) were chosen as follows \cite{Zhang2013511}: $a_0=-3.12513824$, $a_1=1.84108651$, $a_2=-0.35706478$, $a_3= 0.10185626$, $a_4=-0.02924772$, $a_5=0.00696837$, $a_6= -0.00102952$.
\subsection{The solution of the SLAEs}
Let us write down the difference problem  (\ref{main_eq_diff}) in the form of a SLAE as
\begin{equation}
\label{main_slae}
\begin{array}{l}
\left(
\begin{array}{llll}
\gamma_1 c^2\mathcal{L}_{x}-\tilde{\eta}^2I&0&0&\beta_1 c^2\mathcal{L}_{x}\\
0&\gamma_2 c^2 \mathcal{L}_{x}-\tilde{\eta}^2I&0&\beta_2 c^2\mathcal{L}_{x}\\
0&0&\gamma_3 c^2\mathcal{L}_{x}-\tilde{\eta}^2I&\beta_3 c^2\mathcal{L}_{x}\\
-251/720\tilde{\eta} I&-251/720\tilde{\eta} I&-251/720\tilde{\eta} I&\left(c/h_z+251/720\tilde{\eta}\right)I
\end{array}
\right)\left(
\begin{array}{l}
\bar{\mathbf{\Psi}}_{k+1}^{m,1}\\
\bar{\mathbf{\Psi}}_{k+1}^{m,2}\\
\bar{\mathbf{\Psi}}_{k+1}^{m,3}\\
\bar{\mathbf{U}}_{k+1}^m
\end{array}
\right) \\ = \left(
\begin{array}{c}
\Phi_2\left(\bar{\mathbf{\Psi}}_{k+1}^{m,1}\right)\\
\Phi_2\left(\bar{\mathbf{\Psi}}_{k+1}^{m,2}\right)\\
\Phi_2\left(\bar{\mathbf{\Psi}}_{k+1}^{m,3}\right)\\
c/hz\bar{\mathbf{U}}_{k}^m+
1/720\left(\sum_{i=-3}^0\alpha_i\left(\tilde{\eta}\mathbf{\Theta}_{k+i}^m+
\Phi_1\left(\bar{\mathbf{\Theta}}_{k+i}^{m}\right)\right)+\alpha_1\Phi_1\left(\bar{\mathbf{\Theta}}_{k+1}^{m}\right)\right)
\end{array}
\right),
\end{array}
\end{equation}
where $\bar{\mathbf{\Theta}}_{k}^{m}=-\bar{\mathbf{U}}_{k}^m+\sum_{s=1}^3\bar{\mathbf{\Psi}}_{k}^{m,s}$ and $I$ is the unit matrix. Employing the Schur complement \cite{Cottle1974}, the mesh functions $\bar{\mathbf{U}}_{k+1}^m$ can be calculated through the solution to the following reduced SLAE
\begin{equation}
\label{main_schur}
\begin{array}{l}\displaystyle
\left[\left(c/h_z+\frac{251}{720}\tilde{\eta}\right)I+\frac{251}{720}\tilde{\eta}\sum_{s=1}^3\beta_s c^2 \mathcal{L}_{x}\left( \gamma_s c^2 \mathcal{L}_{x}-\tilde{\eta}^2I\right)^{-1}\right]\mathbf{\bar{U}}_{k+1}^m =\mathbf{\bar{F}}_u^m+\sum_{s=1}^3M_s^{-1}\bar{\mathbf{F}}^m_{\psi_s},
\end{array}
\end{equation}
where
$$
\begin{array}{ll}
\displaystyle
M_s=\frac{\gamma_s c^2}{\tilde{\eta}^2} \mathcal{L}_{x}-I,\\
\displaystyle
\mathbf{\bar{F}}_u^m=c/hz\bar{\mathbf{U}}_{k}^m+
\frac{1}{720}\left(\sum_{i=-3}^0\alpha_i\left(\tilde{\eta}\bar{\mathbf{\Theta}}_{k+i}^{m}+
\Phi_1\left(\bar{\mathbf{\Theta}}_{k+i}^{m}\right)\right)+
\alpha_1\Phi_1\left(\bar{\mathbf{\Theta}}_{k+1}^{m}\right)\right),\\
\displaystyle \bar{\mathbf{F}}^m_{\psi_s}=\Phi_2\left(\bar{\mathbf{\Psi}}_{k+1}^{m,s}\right)/\tilde{\eta}^{2}.
\end{array}
$$
Making use of the matrix property \cite{Henderson1981} for (\ref{main_schur})
\begin{equation}
\label{matrix:prop1}
\left(B+I\right)^{-1}B=I-\left(B+I\right)^{-1},
\end{equation}
multiplying the equation by the matrix $M_1M_2M_3$ and taking into consideration the commutative property of \mbox{$M_iM_j=M_jM_i$}, we obtain the governing equation for the calculation of the mesh functions $\mathbf{\bar{U}}^m_{k+1}$
\begin{equation}
\label{reduce_eq}
\begin{array}{ll}
\displaystyle
\left[M_1M_2M_3\left( {c}/{h_z}+\tilde{\eta}+\frac{251}{720}\tilde{\eta}\sum_{s=1}^3\frac{\beta_s}{\gamma_s}\right)I+
\frac{251}{720}\tilde{\eta}\left(\frac{\beta_1}{\gamma_1}M_2M_3+\frac{\beta_2}{\gamma_2}M_1M_3+\frac{\beta_3}{\gamma_3}M_1M_2\right)\right]\mathbf{\bar{U}}^m_{k+1}\\\\
\displaystyle =M_1M_2M_3\mathbf{\bar{F}}_u^m+\tilde{\eta}\left(M_2M_3\bar{\mathbf{F}}^m_{\psi_1}+M_1M_3\bar{\mathbf{F}}^m_{\psi_2}+M_1M_2\bar{\mathbf{F}}^m_{\psi_3}\right).
\end{array}
\end{equation}
As opposed to the Fourier transform, the coefficients of the Laguerre expansion in series (\ref{series_lag}) are dependent in a recurrent manner (\ref{main_eq_laguerre}). Hence, for a fixed $k$ for different $m$ it is required to solve SLAEs many times with a common real matrix and different right-hand sides.
Matrix (\ref{reduce_eq}) is banded and can be explicitly represented without calculation of the matrices $M_s^{-1}$ thus allowing us to apply efficient algorithms for solving SLAEs based on $LU$-decomposition.
Solving the SLAEs with banded matrices with a parallel algorithm, it is reasonable to use the parallel dichotomy algorithm \cite{terekhov:Dichotomy,Terekhov2016,Terekhov:2013}, which was developed for tridiagonal matrices and block-tridiagonal matrices. With respect to the number of arithmetical operations, the dichotomy algorithm is comparable with other available algorithms; however the time needed for inter-process communications is considerably less in the dichotomy algorithm as compared to other algorithms. This is because the implementation of the dichotomy process on a supercomputer reduces to calculating the sum of series for distributed data. The commutative and associative properties of addition enable a considerable reduction in the total computation time with the use of inter-processor interaction optimization algorithms.

After the calculation of the mesh functions $\bar{\mathbf{\mathbf{U}}}^{m}_{k+1}$ before turning to calculating the functions $\bar{\mathbf{\mathbf{U}}}^{m}_{k+2}$, the functions $\bar{\Psi}^{m,s}_{k+1}$ should  be calculated as
\begin{equation}\displaystyle
M_s \bar{\mathbf{\Psi}}^{m,s}_{k+1}=\tilde{\eta}^{-2}\left(-\beta_sc^2\mathcal{L}_{x}\bar{\mathbf{U}}^m_{k+1}+\Phi_2\left(\bar{\mathbf{\Psi}}^{m,s}_{k+1}\right)\right),\; s=1,2,3.
\label{psi_non_reduce}
\end{equation}
Making use of the property (\ref{matrix:prop1}), we arrive at
\begin{equation}
\label{psi_func}
\bar{\mathbf{\Psi}}^{m,s}_{k}=M_s^{-1}\left(-\frac{\beta_s}{\gamma_s}\bar{\mathbf{U}}^m_{k}+\frac{1}{\tilde{\eta}^2}\Phi_2\left(\bar{\mathbf{\Psi}}^{m,s}_{k}\right)\right)
-\frac{\beta_s}{\gamma_s}\bar{\mathbf{U}}^m_{k},\quad s=1,2,3.
\end{equation}
For solving the 2D  OWWE one needs not only the stabilization of the numerical instability of difference approximation for the operator $\partial/\partial z$, but also the instability of the real Pad\'e approximation (\ref{sqr-root-approx}). The method AM5-I5 allows solving both these problems.

{\bf Algorithm 2-1. The Adams-Moulton downward-continuation procedure for the 2D OWWE.}

To calculate the mesh functions $\bar{\mathbf{U}}^m,\bar{\mathbf{\Psi}}^m$ accurate to $O(h_x^{\xi}+h_z^5)$, the following is necessary:
 \begin{enumerate}
 \item    Let the number of nodes of the mesh $\Omega$ in the direction $z$ be odd.  For all $i$ for the functions
$\Phi_1(\bar{u}^m_{ik})$, $\Phi_2(\bar{\psi}^{m,s}_{ik})$ construct, independently, the 1D quintic splines in the direction $z$, using only odd values of   $k$.

  \item   Replace the values of the functions $\Phi_1(\bar{u}^m_{ik})$, $\Phi_2(\bar{\psi}^{m,s}_{ik})$ for even $k$ by their interpolated values (the quintic spline filtration).
   \item For $k=4,...,K-1$
   \begin{enumerate}
   \item [3.1.] Applying equation (\ref{reduce_eq}), calculate the solution $\bar{\mathbf{U}}^{m}_{k+1}$.
   \item [3.2.] Applying equation (\ref{psi_non_reduce}), calculate the solution $\bar{\mathbf{\Psi}}^{m,s}_{k+1},\;s=1,2,3$
   \end{enumerate}
   \item Turn to the calculation of the  $(m+1)\text{th},\;\text{the} (m+2)\text{th}$, etc. coefficients of the expansion of the Laguerre series.
\end{enumerate}

The above-considered way of stabilizing the solution allows the stability not only of  the Adams-Moulton implicit schemes, but also the Adams-Bashfort explicit schemes that are stable for essentially lesser steps $h_z$. As consequence, the number of SLAEs to be solved of the form (\ref{reduce_eq}),(\ref{psi_non_reduce}) is multiply increased thus making the Adams-Bashfort method inefficient. Let us now consider the predictor-corrector method combining the computational efficiency of explicit and high stability of the implicit Adams schemes.

For equation (\ref{main_eq1_laguerre}) as a predictor procedure we choose the Adams-Bashfort explicit method of the fifth order of accuracy
\begin{equation}
\frac{\bar{u}^m_{ik+1}-\bar{u}^m_{ik}}{h_z}=\frac{1}{720c}\sum_{j=-4}^0\varrho_j\left(\sum_{s=1}^{3}\left(\tilde{\eta}\bar{\psi}^{m,s}_{ik+j}
+\Phi_1\left(\bar{\psi}^{m,s}_{ik+j}\right)\right)-\tilde{\eta}\bar{u}^m_{ik+j}-\Phi_1\left(\bar{u}^m_{ik+j}\right)\right),
\label{Adams_Bashfort}
\end{equation}
where  $\varrho_0=1901,\;\varrho_{-1}=-2774,\;\varrho_{-2}=2616,\;\varrho_{-3}=-1274,\;\varrho_{-4}=251$.
In terms of correction we will use the Adams-Moulton scheme of the fifth order of accuracy (\ref{adams_owwe.a}), where unknown values $\bar{\psi}^{m,s}_{ik+1}$  should be replaced by predicted values. In this case, the functions $\bar{u}^m_{ik+1}$ are explicitly expressed. If for correction we use scheme (\ref{adams_owwe.a})  and substitute the predicted values both for $\bar{\psi}^{m,s}_{ik+1}$  , and for $\bar{u}^m_{ik+1}$ into the right-hand side, then for the sake of stability, smaller steps $h_z$  and simultaneously a larger number of corrections will be required.

{\bf Algorithm 2-2. The Predictor-Corrector downward-continuation procedure for the 2D OWWE.}

To calculate the mesh functions $\bar{\mathbf{U}}^m,\bar{\mathbf{\Psi}}^m$ accurate to $O(h_x^{\xi}+h_z^5)$, the following is necessary:
 \begin{enumerate}
   \item  Let the number of nodes of the mesh $\Omega$ in the direction $z$ be odd.  For all $i$ for the functions $\Phi_1(\bar{u}^m_{ik})$, $\Phi_2(\bar{\psi}^{m,s}_{ik})$ independently construct the 1D quintic splines in the direction  $z$ using only odd values of $k$.
  \item   Replace values of the functions $\Phi_1(\bar{u}^m_{ik})$, $\Phi_2(\bar{\psi}^{m,s}_{ik})$ for even $k$ by their interpolated values  (the quintic spline filtration).
   \item For  $k=4,...,K-1$
   \begin{enumerate}
   \item [3.1.] Applying equation (\ref{Adams_Bashfort}), calculate the predicted solution ${\breve{\mathbf{U}}}^{m}_{k+1}$.
   \item [3.2.] Applying equation (\ref{psi_func}) with  ${\breve{\mathbf{U}}}^{m}_{k+1}$, calculate the predicted solution $\breve{\Psi}^{m,s}_{k+1},\;s=1,2,3$ .
   \item [3.3.] Applying equation (\ref{adams_owwe.a}), substituting  $\breve{\Psi}^{m,s}_{k+1}$ instead of $\bar{\Psi}^{m,s}_{k+1}$, calculate the corrected solution $\breve{\breve{\mathbf{U}}}^{m}_{k+1}$.
   \item [3.4.] Applying equation (\ref{psi_func}) with $\breve{\breve{\mathbf{U}}}^{m}_{k+1}$, calculate the final solution for $\bar{\mathbf{\Psi}}^{m,s}_{k+1},\;s=1,2,3$ .
   \item [3.5.] Applying equation (\ref{adams_owwe.a}) with $\bar{\mathbf\Psi}^{m,s}_{k+1}$, calculate the final solution for $\bar{{\mathbf{U}}}^{m}_{k+1}$.
   \end{enumerate}
   \item Turn to the calculation of the $(m+1)\text{th},\;\text{the }(m+2)\text{th}$, etc. coefficients of the expansion of the Laguerre series.
\end{enumerate}

Thus, within the predictor-corrector method instead of indefinite non-symmetrical SLAE  (\ref{reduce_eq}) it is necessary to solve SLAEs of the form (\ref{psi_func}) with sign-defined symmetric matrices of lesser dimensions that allows the use of efficient algorithms of the computational linear algebra and a decrease in the calculation time.  As compared to the Marchuk-Strang method that is accurate to the second order, the predictor-corrector method is of the fifth order of accuracy. Further increase of approximation order is not reasonable because for providing the stability an essentially smaller step  $h_z$ should be set.

Both for the Adams method and for the predictor-corrector method one needs initial values to start the calculation. It is required to use other methods such as the Richardson extrapolation or Crank-Nicolson scheme with a smaller step for obtaining such  initial values.
\section{Numerical experiments for the 2D one-way wave equation}
Analytically it is really difficult to provide a strict mathematical substantiation of algorithms of  a high-order of accuracy for the
2D OWWE. Therefore to confirm the efficiency of the approaches proposed,  thorough testing is needed.
Let us discuss a few tests that would allow the evaluation of
the quality of the solution to be obtained as compared with the Richardson algorithm. Numerical procedures were implemented in Fortran-90 using the MPI library.

\subsection{The impulse response}
In the first test we illustrate analyzing the accuracy by the impulse responses. For the calculation, we used the homogeneous medium model with the speed $250\;m/s$ and the size $3.5\;km\times1.5\;km$. The point source (\ref{source}) with the parameters $t_0=0.2\,s,\,\delta=4,\,f_0=30\,\text{Hz}$ was located at the center of the upper surface.
The number of addends in series (\ref{series_lag}) was $n=4000$ for $T=6\;s$; the expansion parameter was $\eta=600$.

A disadvantage of the Laguerre transform is the absence of the fast transformation algorithm. The implementation of the forward Laguerre transform (\ref{series_lag}) with the help of the method of least squares, the number of arithmetical operations will be of order $O(NP)$, where $P$ is the number of discrete points of the approximated function and $N\leq P$ is the number of terms in the Laguerre series needed for attaining the required accuracy in the norm $L_2$, whereas for the fast Fourier transform computer costs are essentially less and make up $O(N\log(N))$. However, taking into account the fact that input data are set only along the upper surface ($z=0$), and the inverse transformation is done for a fixed time instant, the total cost of the direct and inverse transformations  appears to be minor as compared to that needed for calculation of coefficients of series (\ref{series_lag}) from the solution to problem (\ref{main_eq_laguerre}). The numerical experiments have confirmed that the time needed for carrying out the Laguerre transform for the initial data is less than one percent of the total calculation time.

\begin{figure}[!htb]
\centering
        \begin{subfigure}[b]{0.49\textwidth}
                \includegraphics[width=\textwidth]{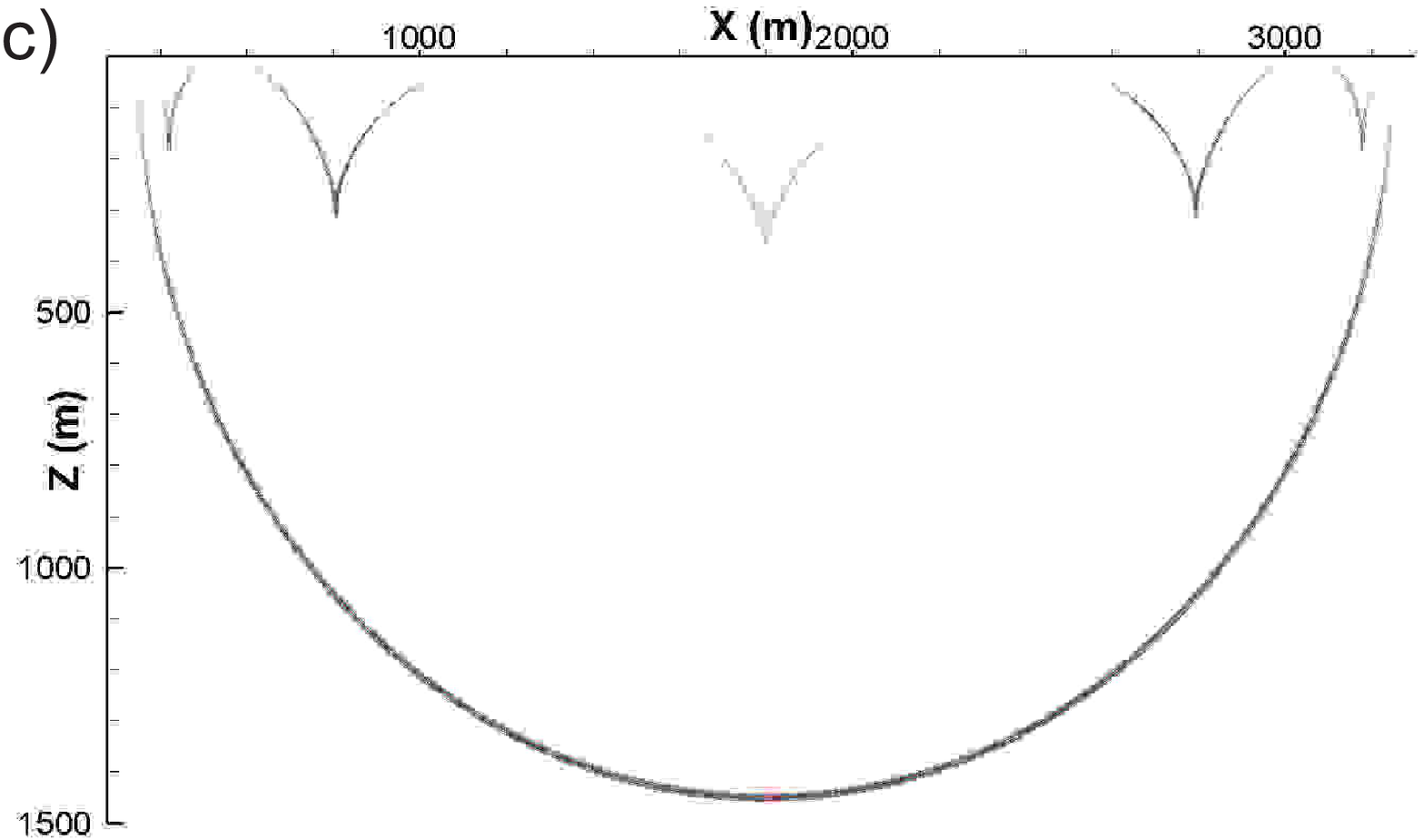}
                \label{fig:Test1-One-Way-Rich}
        \end{subfigure}
        ~
        \begin{subfigure}[b]{0.49\textwidth}
                \includegraphics[width=\textwidth]{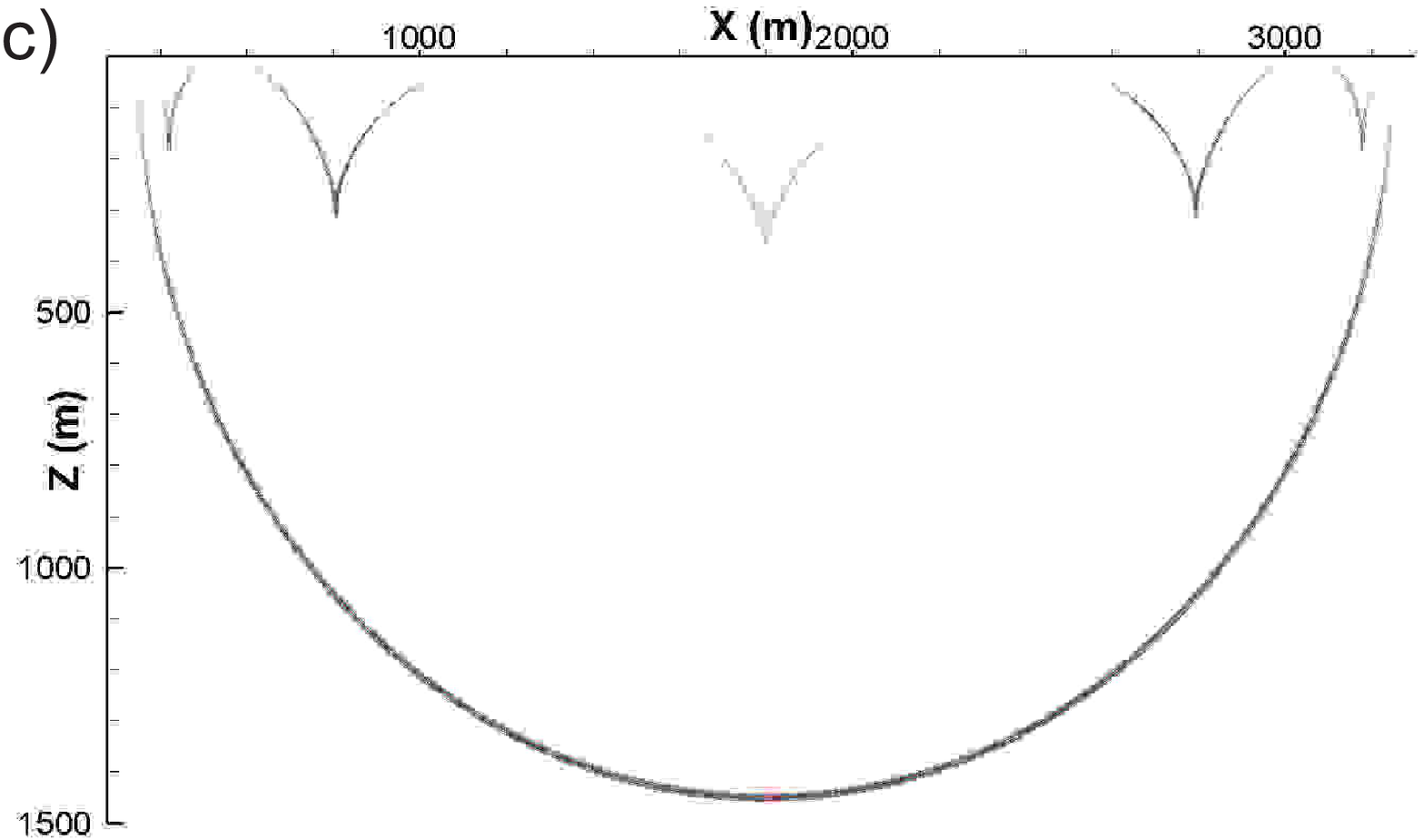}
                \label{Test1-One-Way-Rich2}
        \end{subfigure}

        \begin{subfigure}[b]{0.49\textwidth}
                \includegraphics[width=\textwidth]{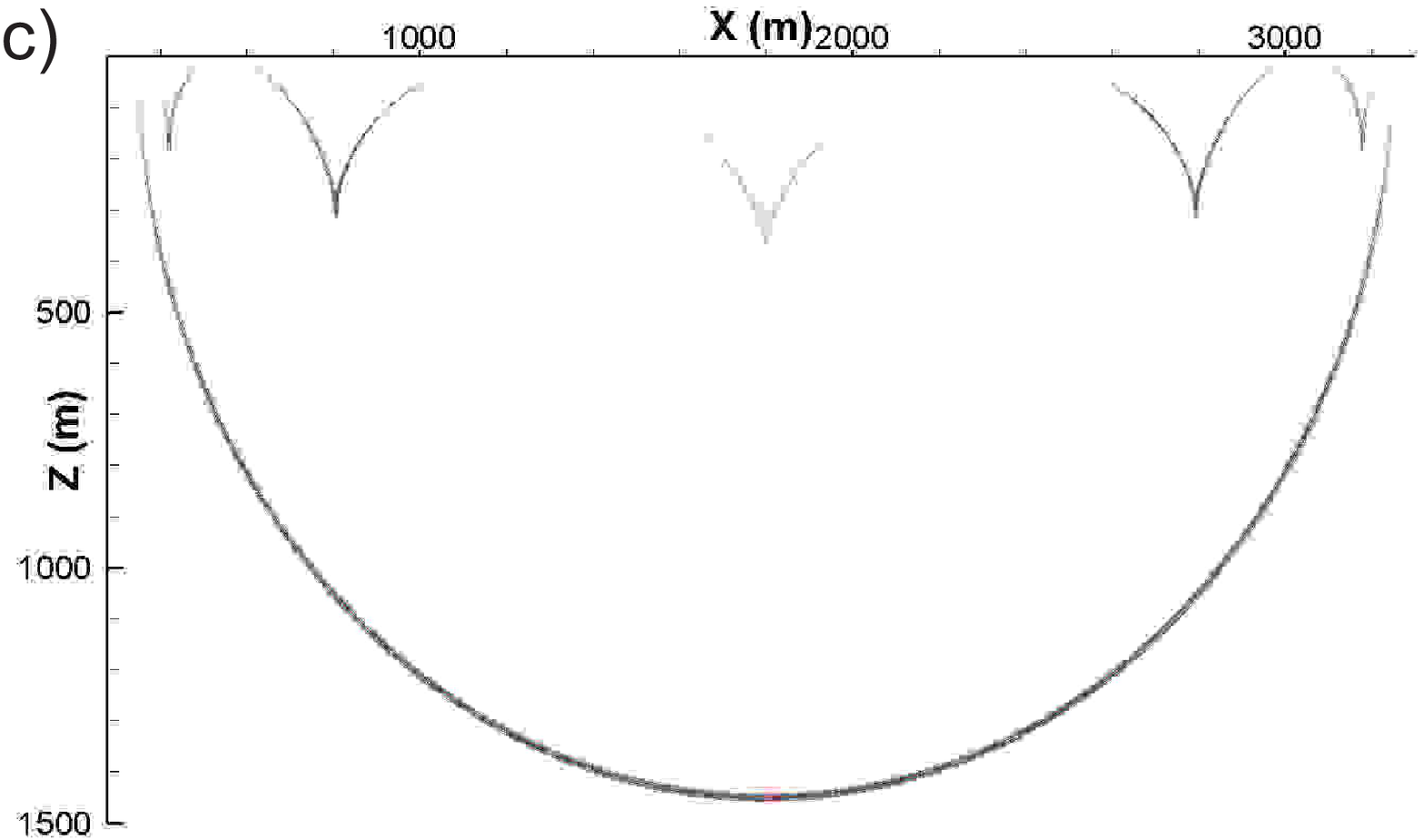}
                \label{Test1-One-Way-Adams}
        \end{subfigure}
        ~
        \begin{subfigure}[b]{0.49\textwidth}
                \includegraphics[width=\textwidth]{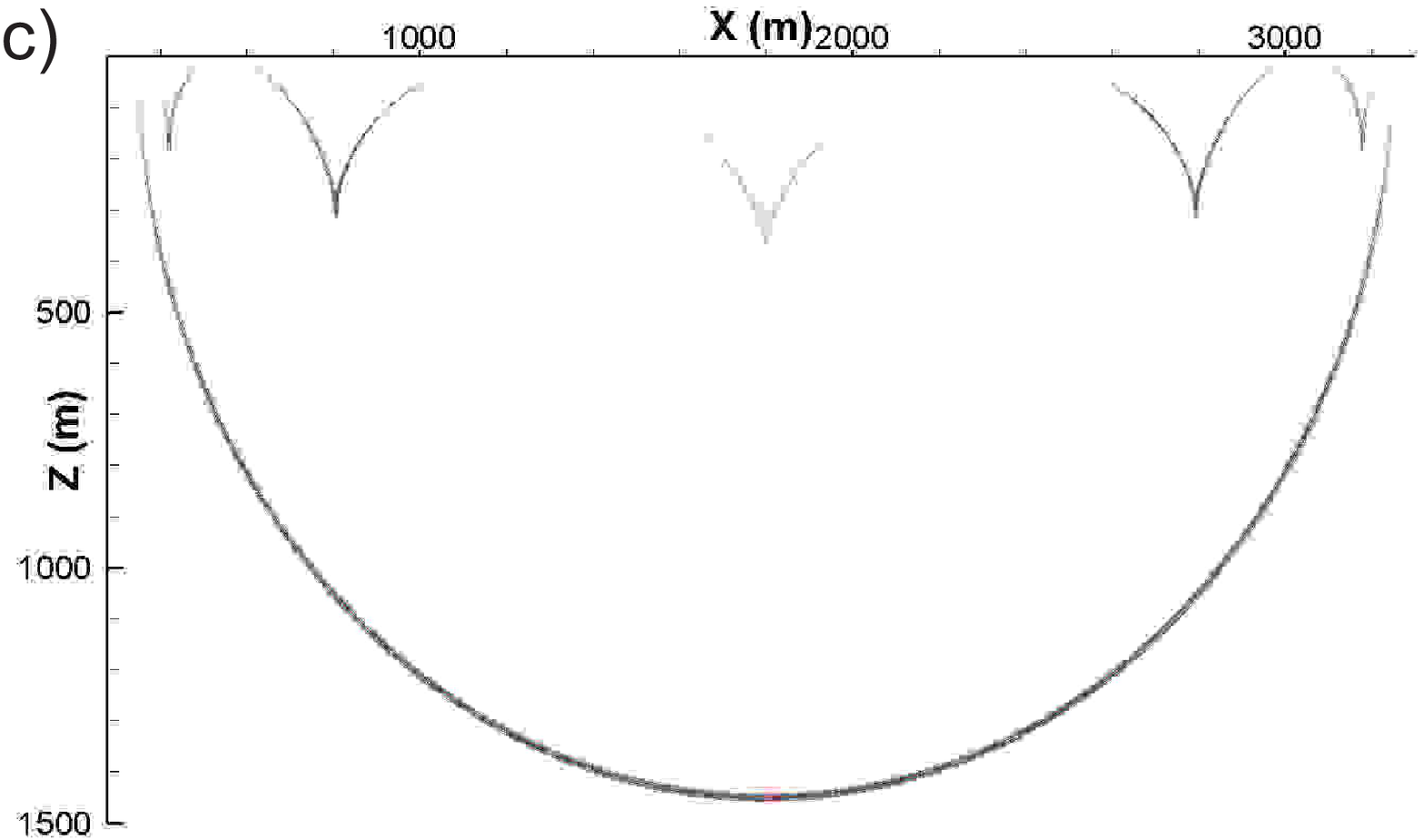}
                \label{Test1-One-Way-Predictor}
        \end{subfigure}

\caption{Snapshots for the wave field at $t=6\;s$ for the homogeneous velocity model. The Richardson extrapolation (a)~\mbox{$h_{x,z}=1\;m$} and (b)~\mbox{$h_{x,z}=0.5\;m$}, (c)~AM5-I5 method with \mbox{$h_{x,z}=1\;m$}, (e)~PC5-I5 method with \mbox{$h_{x,z}=1\;m$}.}
\label{fig:impulse-one-way}
\end{figure}

\begin{figure}[!htb]
\centering
\includegraphics[width=0.49\textwidth]{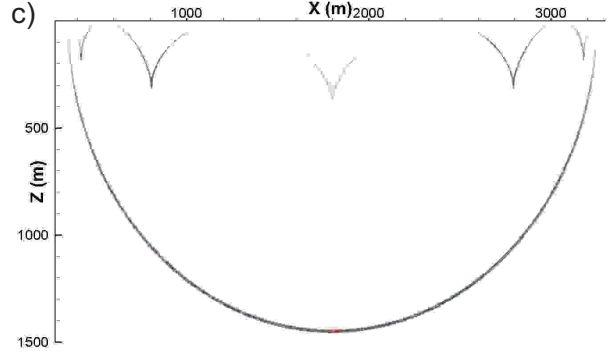}
\caption{ Dependence of the wave field on the coordinate along the straight line "Slice"\;(Fig.~\ref{fig:impulse-one-way}) for different meshes and methods. }
\label{Test1-One-Way-Compare}
\end{figure}

In Section \ref{sec:1dowwe}, for the 1D OWWE it was shown that with an equal mesh size the Richardson extrapolation is more accurate than the  AM5-I5-scheme, which is also valid for the 2D case. If one selects  the step $h_z=1$ m for the mesh $\Omega_1$ for the Richardson method and the step $h_z=0.3$ m for the AM5-15 and PC5-15 methods, then the whole volume of calculations and calculation accuracy for all the three methods will be compatible (Fig.~\ref{fig:impulse-one-way}). However with a thorough consideration of values of the amplitudes along the straight line "Slice", it appears to be clear (Fig.~\ref{Test1-One-Way-Compare}) that the AM5-I5 and the PC5-I5 methods are more accurate when the number of nodes of the mesh $\Omega$ is equal to the number of nodes of the meshes $\Omega_1,\Omega_2$.

A considerable accuracy and computational efficiency of  the  AM5-I5, PC5-I5 methods is attained at the expense of lesser stability as compared to the Richardson method. It was experimentally revealed that for the Richardson method the condition of stability is of the form $h_z/h_x\le 1$, while for the  AM5-I5 method it should be $h_z/h_x<0.4$, and for the  PC5-I5 algorithm the stability is attained at  $h_z/h_x<0.3$. The fact that minimum steps required for providing a good accuracy and stability almost coincide, allows us to propose the statement  about the balance of  the PC5-I5 method which is by $30-40\%$ more efficient as compared to the  AM5-I5 method.

Additionally we have considered the AM6-I7 and the AM5-D4 methods, which are stable for the 1D OWWE but unstable for the 2D OWWE. This is because in addition to the numerical instability due to the choice of approximation for the operator $\partial /\partial z$, there is instability caused by the presence of a singular component in the solution to the OWWE, i.e. when denominators (\ref{sqr-root-approx})  are close to zero or vanish.  In the AM5-D4 method, to approximate $\partial /\partial z$ different difference schemes were used for the right-hand and the left-hand sides of equations (\ref{advection_laguerre}),(\ref{advection_laguerre2}) thus additing numerical dissipation, but this approach does not allow restricting the growth of the number of singular components for the  2D OWWE. Also, the AM6-I7 method demonstrates (Fig.~\ref{fig_u2}) lesser dissipation as compared to the algorithm AM5-I5 and, hence,  insufficient level of fictitious absorption does not allow stabilizing the instability for the  2D OWWE.
\subsection{Migration procedures}
 Testing the algorithms based on the combination of the Laguerre transform and the Richardson extrapolation was done in  \cite{Terekhov2017}, where on solving a problem of seismic prospecting, as an example, the post-stack migration procedure within the "explosive boundaries" model has been implemented [4]. As compared to the Finite Difference (FD) \cite{Claerbout:1985}, Fourier Finite Difference (FFD) \cite{Ristow01121994} and Phase Shift Plus Interpolation (PSPI) \cite{Gazdag1984}  methods the algorithm based on the Laguerre transform made possible to obtain a more qualitative solution. Let us now consider similar tests for the algorithms  AM5-I5, PC5-I5 for solving  the 2D OWWE.
\subsubsection{The syncline model}
\begin{figure}[!h]
\centering
        \begin{subfigure}[b]{0.49\textwidth}
                \includegraphics[width=\textwidth]{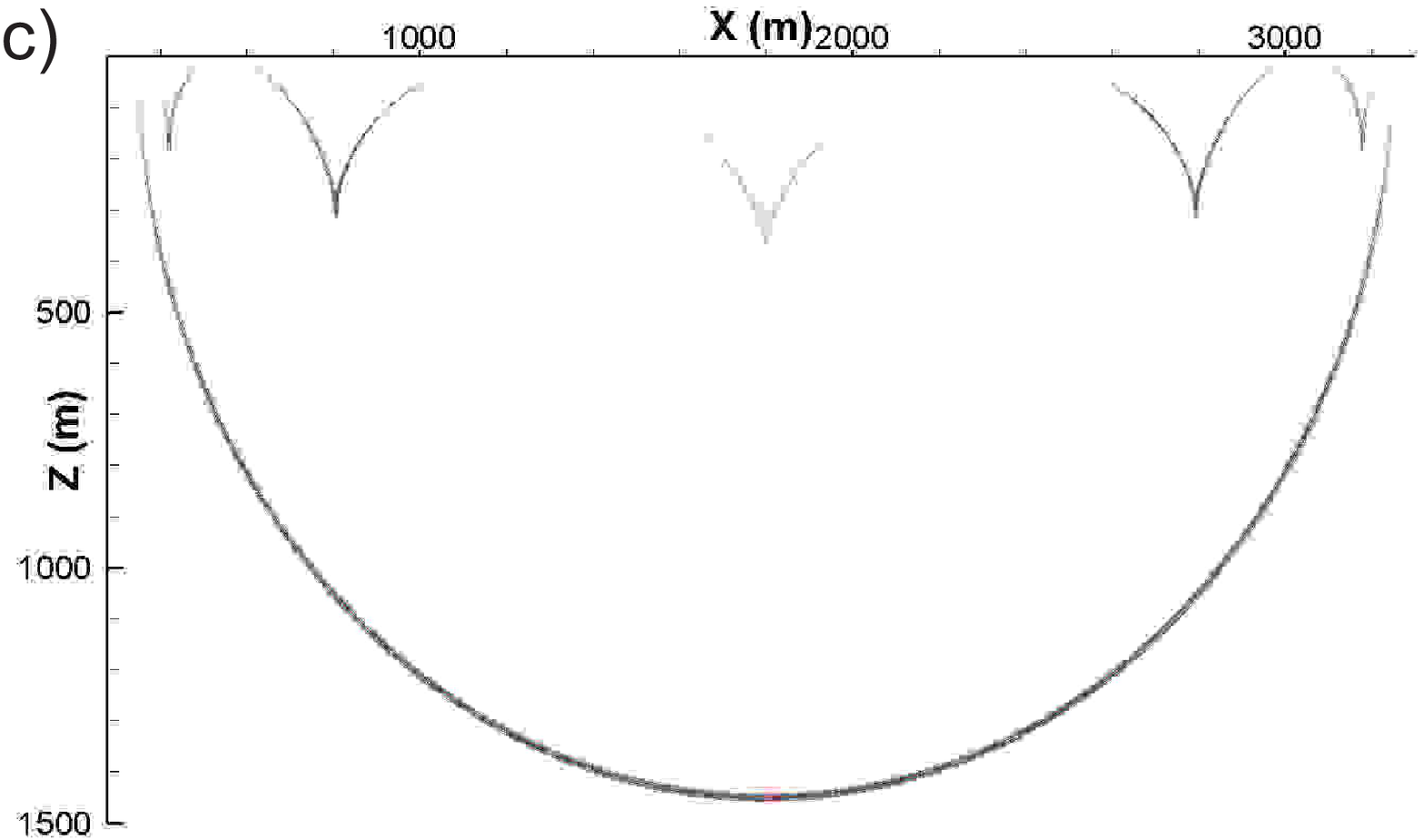}
        \end{subfigure}
              ~
                \begin{subfigure}[b]{0.49\textwidth}
                \includegraphics[width=\textwidth]{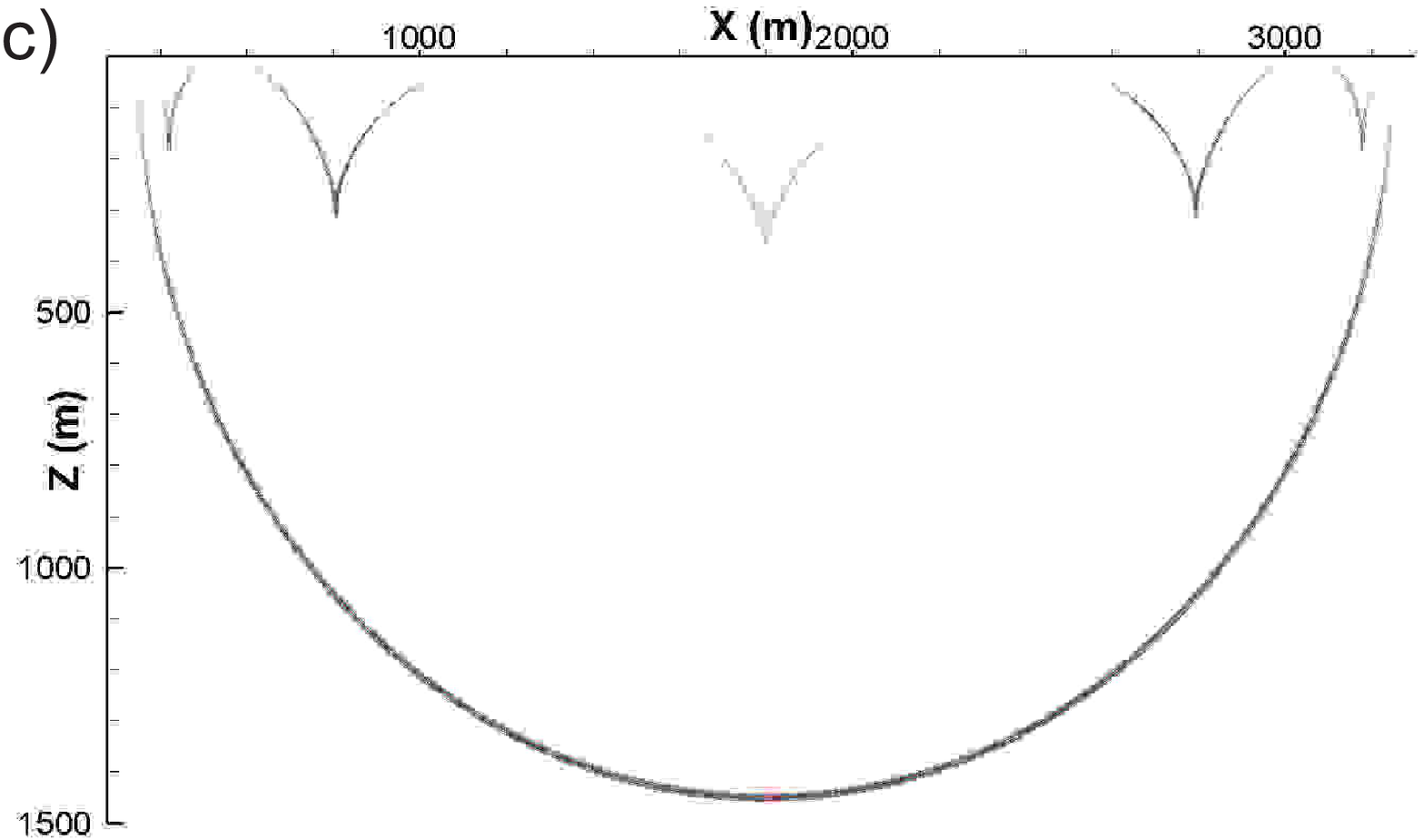}
        \end{subfigure}%
        \caption{(a) Syncline model and (b) zero-offset section.}
        \label{fig:syncline}
\end{figure}

Theoretical seismograms (Fig.~\ref{fig:syncline}b) for the syncline model (Fig.~\ref{fig:syncline}a) were obtained with the help of the Gaussian beams algorithm \cite{GaussianBeams_popov,Cerveny1985} implemented in the package Seismic Unix. For setting the boundary condition on the upper surface, the function for the zero-offset section $\left.u(x,z,t)\right|_{z=0}=g(x,t)$ was expanded in series (\ref{series_lag}) with the parameters $n=2500$ and $\eta=800$ for $t\in[0,4]\;s$.
The calculations were carried out on meshes with the steps $h_{x}=6\;m$, $h_z=3$ or $6$ m for the Richardson method $h_{x}=6\;m$ and $h_z=3$ or $1$ m for the methods AM5-I5, PC5-I5. According to the model of explosive boundaries, the calculation velocities were set to be half the true velocity of the medium model.

The possibility to steadily calculate the solution to the 2D OWWE  for non-smooth velocity models is the advantage of the Richardson method while preliminary smoothing for providing stability is needed for the algorithms AM5-I5, PC5-I5, because application of methods of a high-order is based on the assumption that the solution and the velocity function possess the required smoothness. In order to restrict an increase in numerical instability with a minimum influence on the wave kinematics we use averaging of the form
\begin{equation}
\tilde{c}(i,j)=\frac{1}{8}\left(4c(i,j)+c(i+1,j)+c(i-1,j)+c(i,j+1)+c(i,j-1)\right).
\label{smooth}
\end{equation}
However with decreasing the step $h_x$ the growth of singular  components of the field in the area of  discontinuities of the velocity function increases, therefore a multiple application of formula  (\ref{smooth}) may be needed.
\begin{figure}[!htb]
\centering
        \begin{subfigure}[b]{0.49\textwidth}
                \includegraphics[width=\textwidth]{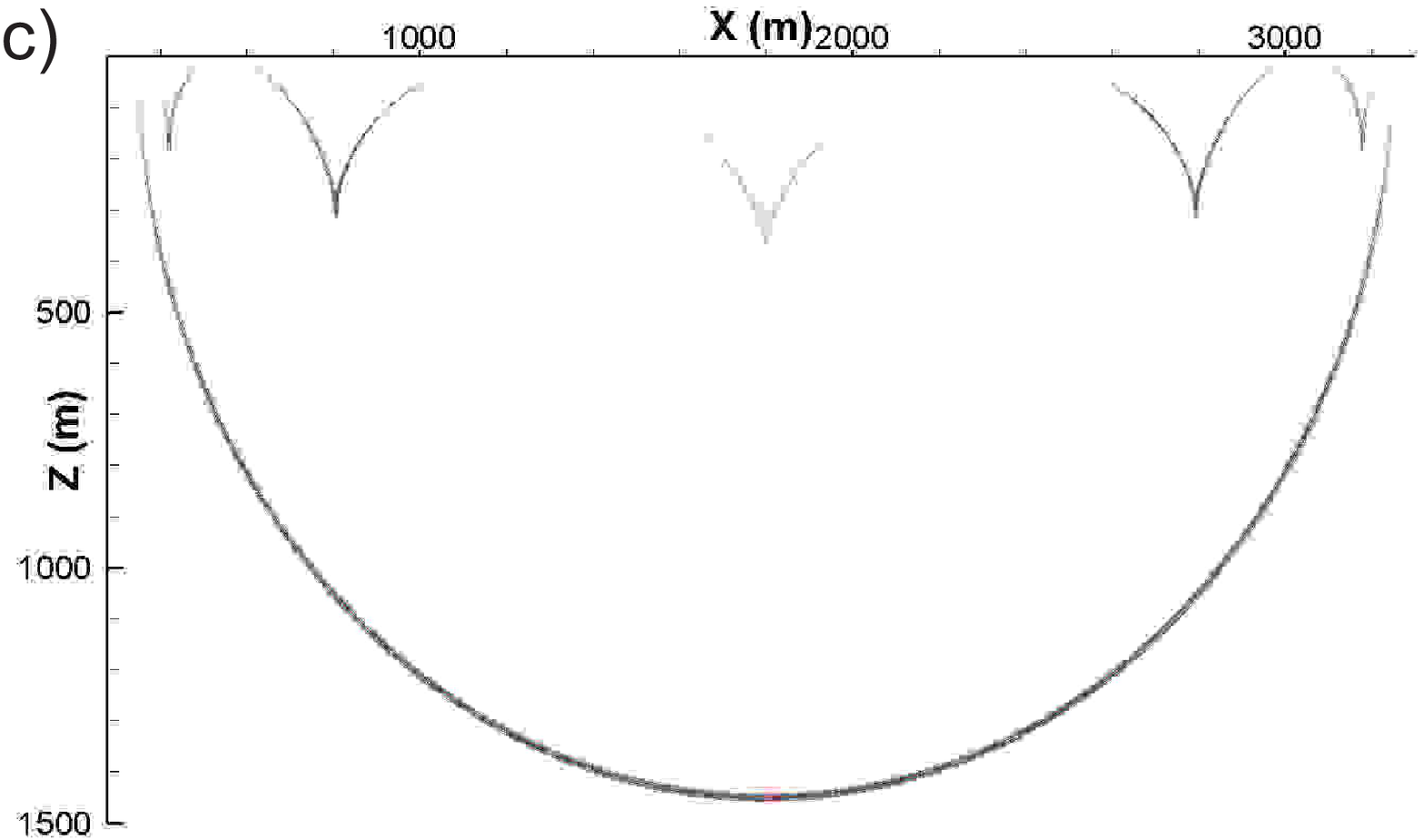}
                \label{fig:Test2-One-Way-Rich}
        \end{subfigure}
        ~
        \begin{subfigure}[b]{0.49\textwidth}
                \includegraphics[width=\textwidth]{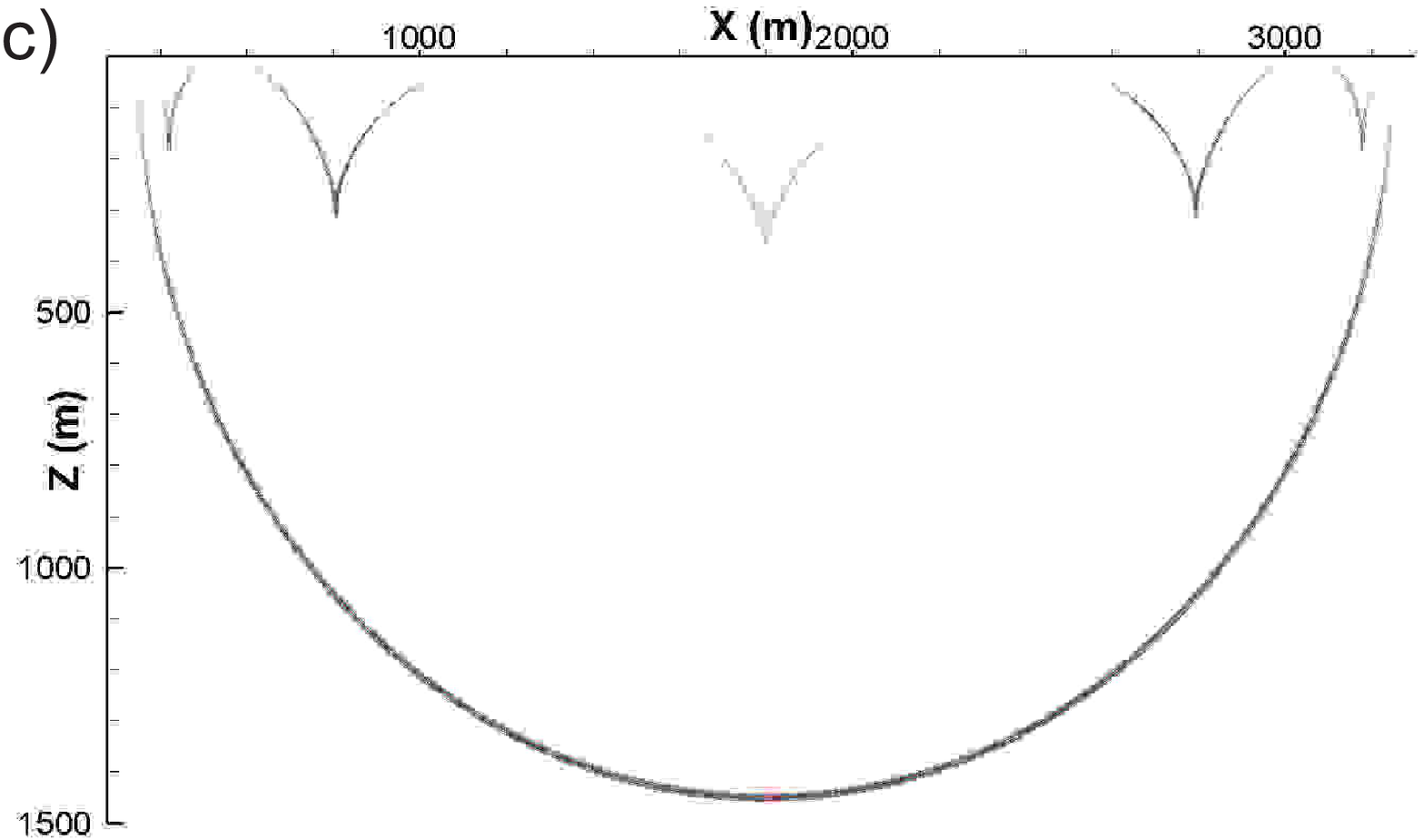}
                \label{Test2-One-Way-Rich2}
        \end{subfigure}

        \begin{subfigure}[b]{0.49\textwidth}
                \includegraphics[width=\textwidth]{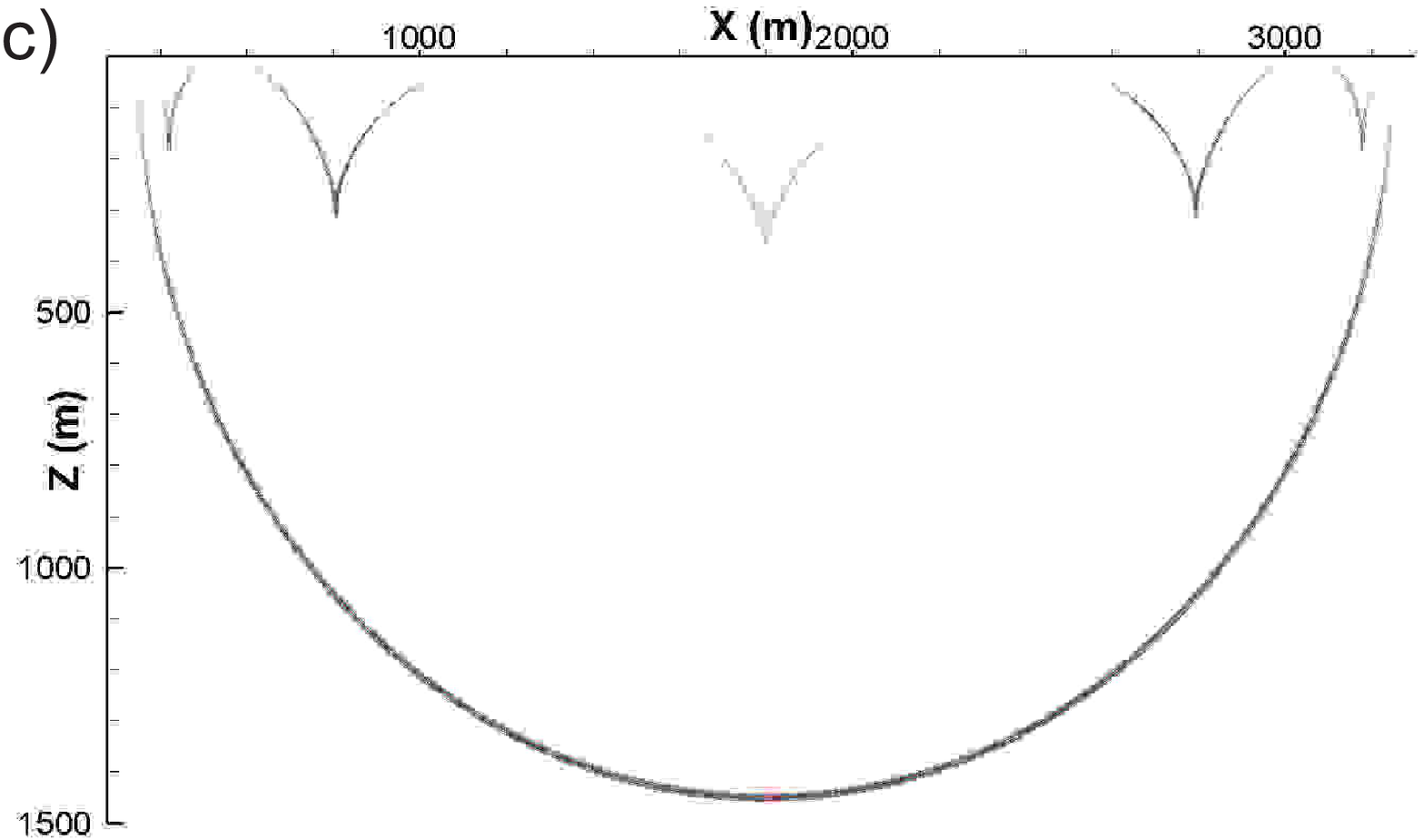}
                \label{Test2-One-Way-Adams}
        \end{subfigure}
        ~
        \begin{subfigure}[b]{0.49\textwidth}
                \includegraphics[width=\textwidth]{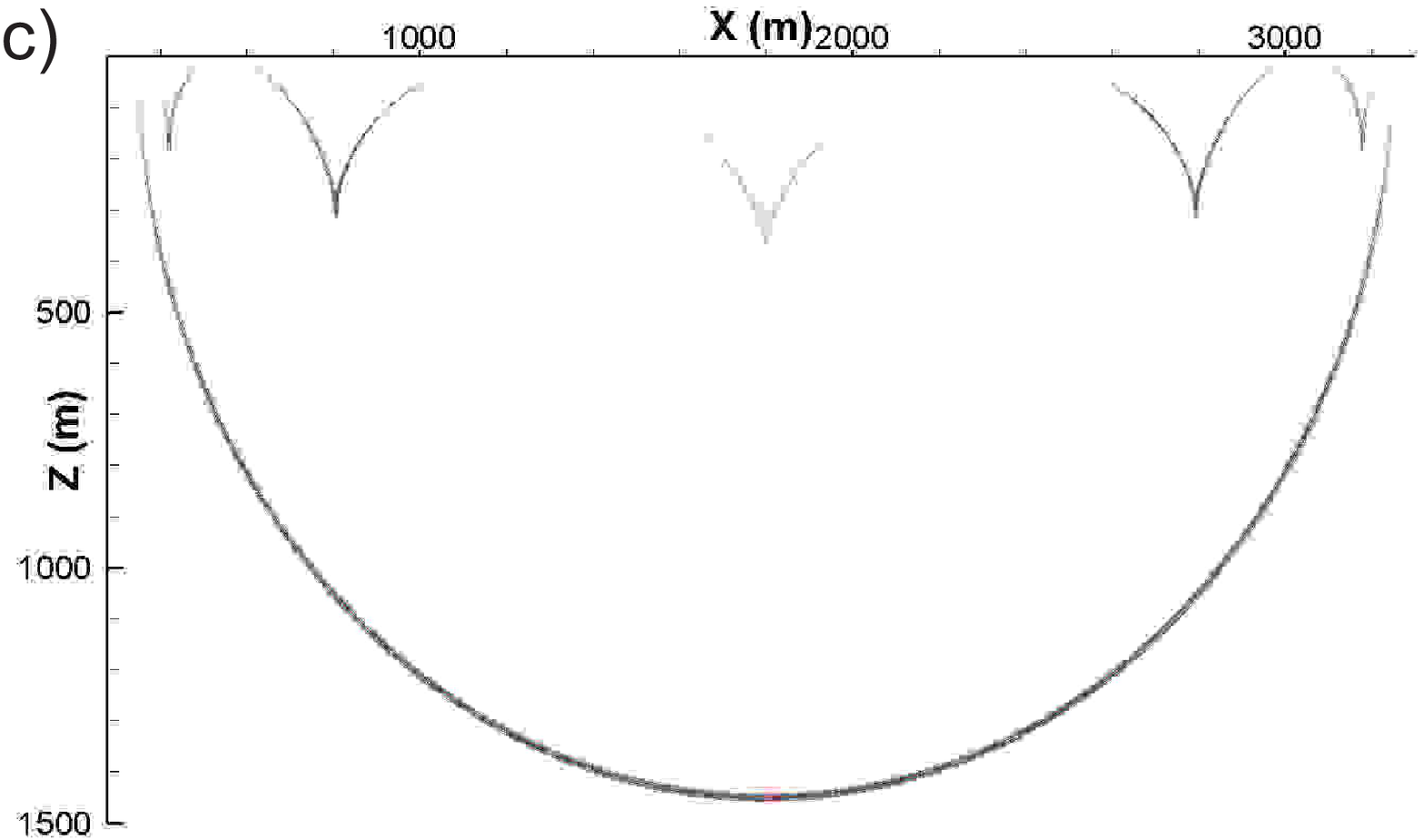}
                \label{Test2-One-Way-Adams2}
        \end{subfigure}

         \begin{subfigure}[b]{0.49\textwidth}
                \includegraphics[width=\textwidth]{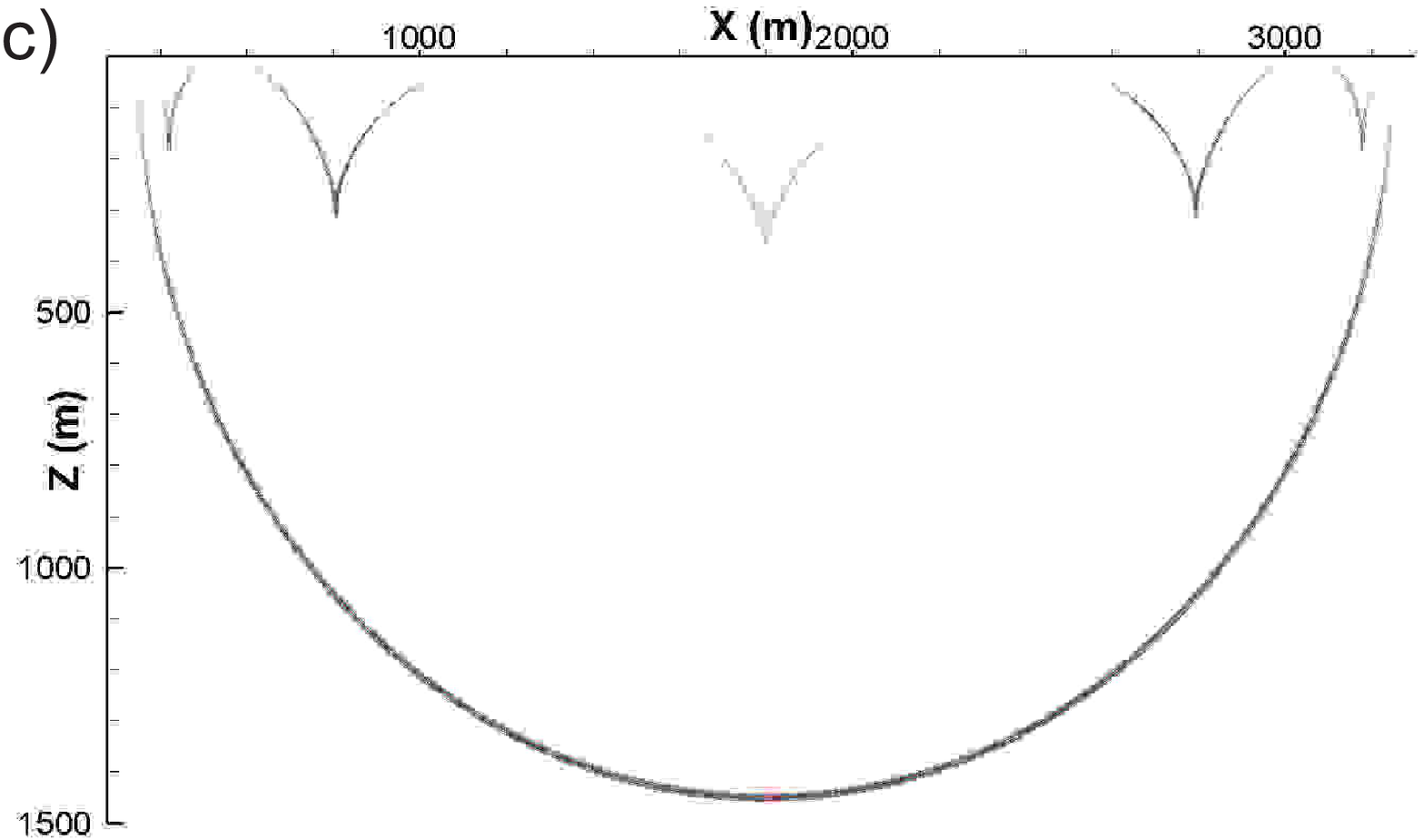}
                \label{Test2-One-Way-PC}
        \end{subfigure}
        ~
        \begin{subfigure}[b]{0.49\textwidth}
                \includegraphics[width=\textwidth]{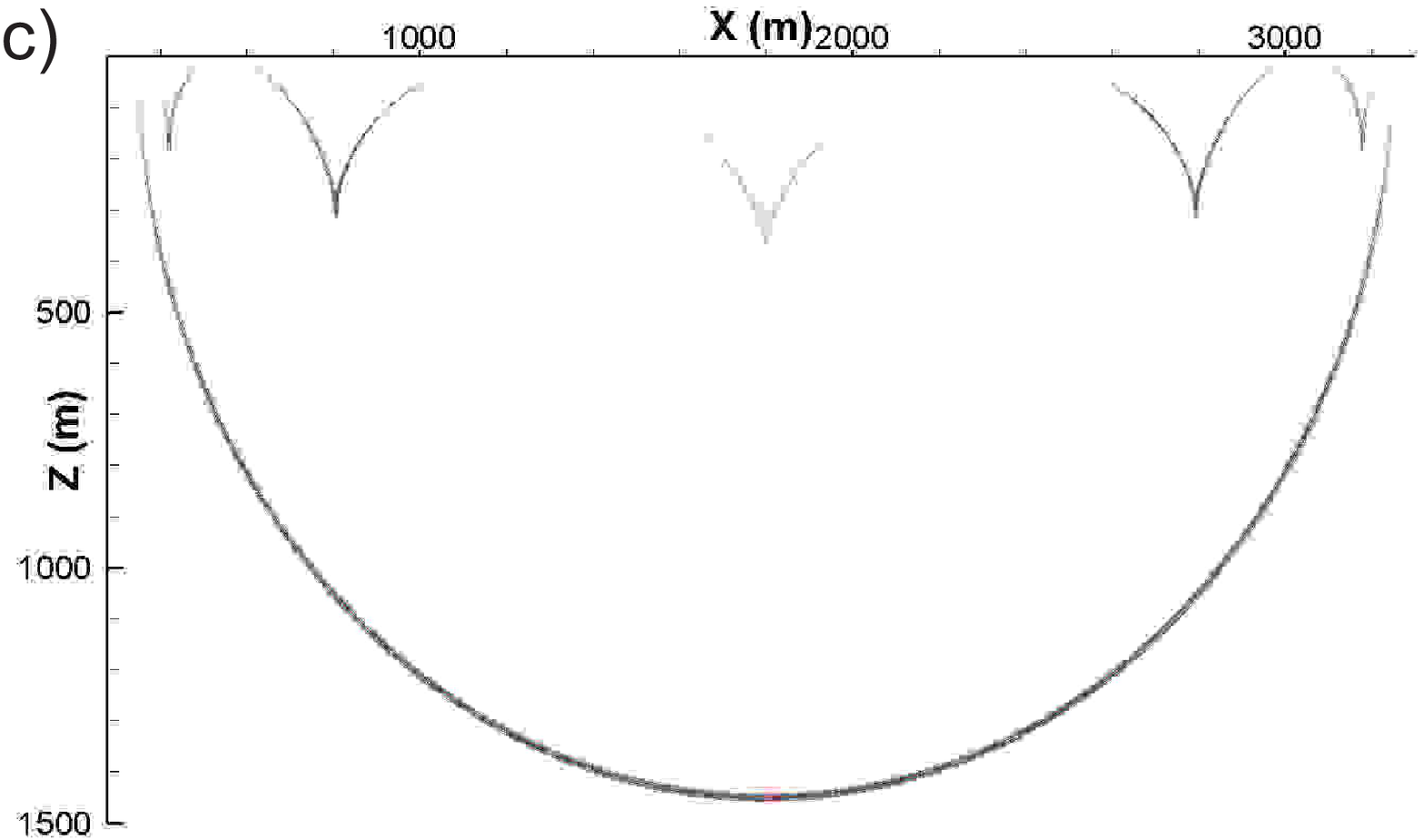}
                \label{Test2-One-Way-PC2}
        \end{subfigure}
\caption{Snapshots for the wave field at $t=4\;s$ for the model and zero-offset section (Fig.~\ref{fig:syncline}) for the Richardson extrapolation (a)~\mbox{$h_{x}=6\;m,\; h_z=6\;m$}, (b)~\mbox{$h_{x}=6\;m, h_z=3\;m$}; for the AM5-I5 method (c)~\mbox{$h_{x}=6\;m,\; h_z=3\;m$}, (d)~\mbox{$h_{x}=6\;m,\; h_z=1\;m$}; for the PC5-I5 method (e)~\mbox{$h_{x}=6\;m,\; h_z=3\;m$}, (f)~\mbox{$h_{x}=6\;m,\; h_z=1\;m$}.   }
\label{fig:syncline_result}
\end{figure}

In Fig.~\ref{fig:syncline_result}a,b it is clear that the Richardson extrapolation is less accurate if the general number of nodes of the meshes  $\Omega_1,\Omega_2$  is equal to the number of nodes of the mesh  $\Omega$ for the methods  AM5-I5, PC5-I5. If the mesh steps are equal, the quality of the image obtained is approximately the same, but the Richardson method requires three times as many calculations as compared to other techniques. The algorithms AM5-I5 (Fig.\ref{fig:syncline_result}c,d) and  PC5-I5(Fig.~\ref{fig:syncline_result}e,f) allow obtaining images that are of the same accuracy, which witnesses to the correctness of calculations because different ideas of using the Adams schemes underlie these methods. The spline-filtration procedure makes possible to provide  the stability with lesser costs than the Richardson extrapolation and the method PC5-I5 is three times more efficient  than the method AM5-I5.
\subsubsection {The Sigsbee model}
\begin{figure}[!h]
\centering
        \begin{subfigure}[b]{0.49\textwidth}
                \includegraphics[width=\textwidth]{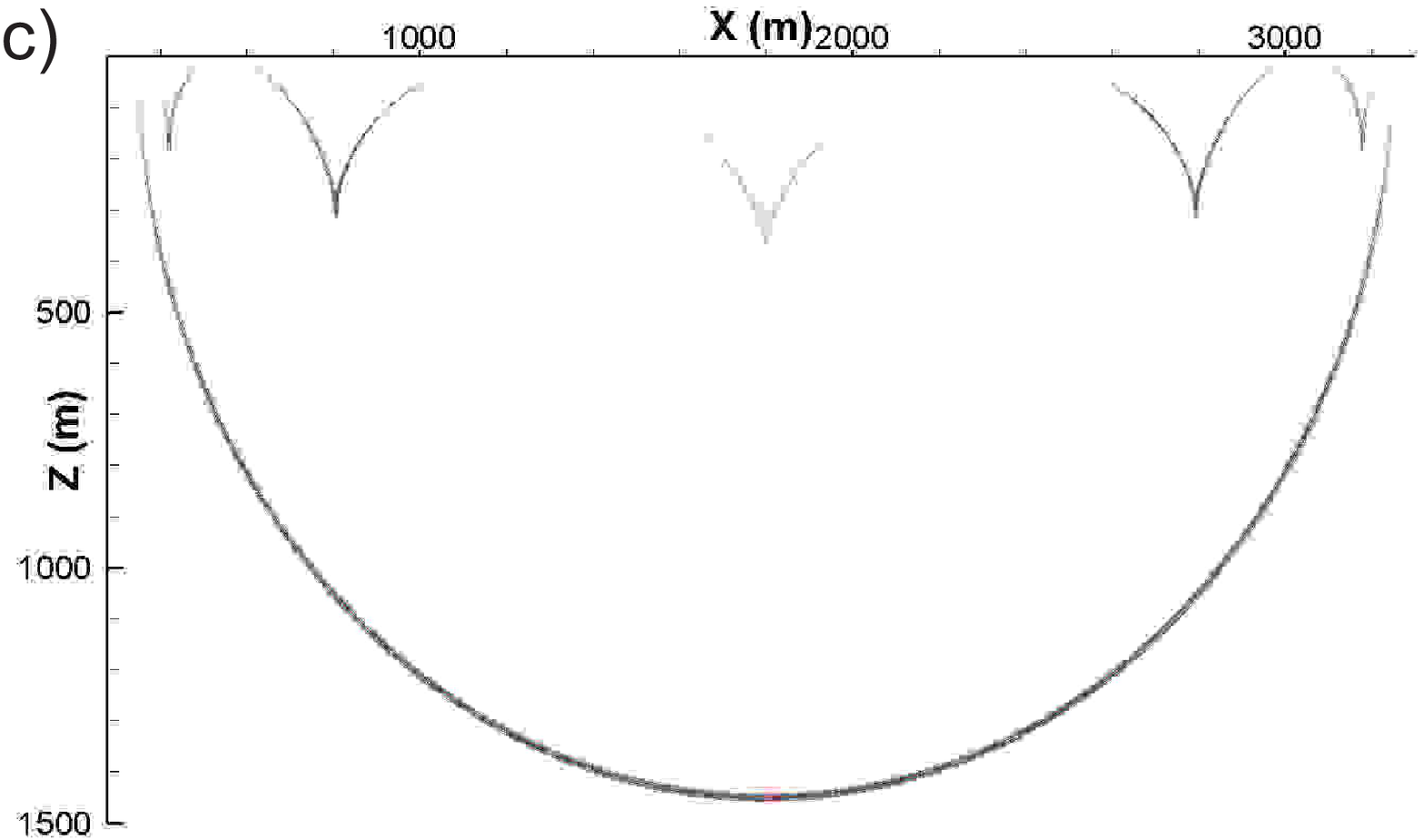}
        \end{subfigure}
              ~
                \begin{subfigure}[b]{0.49\textwidth}
                \includegraphics[width=\textwidth]{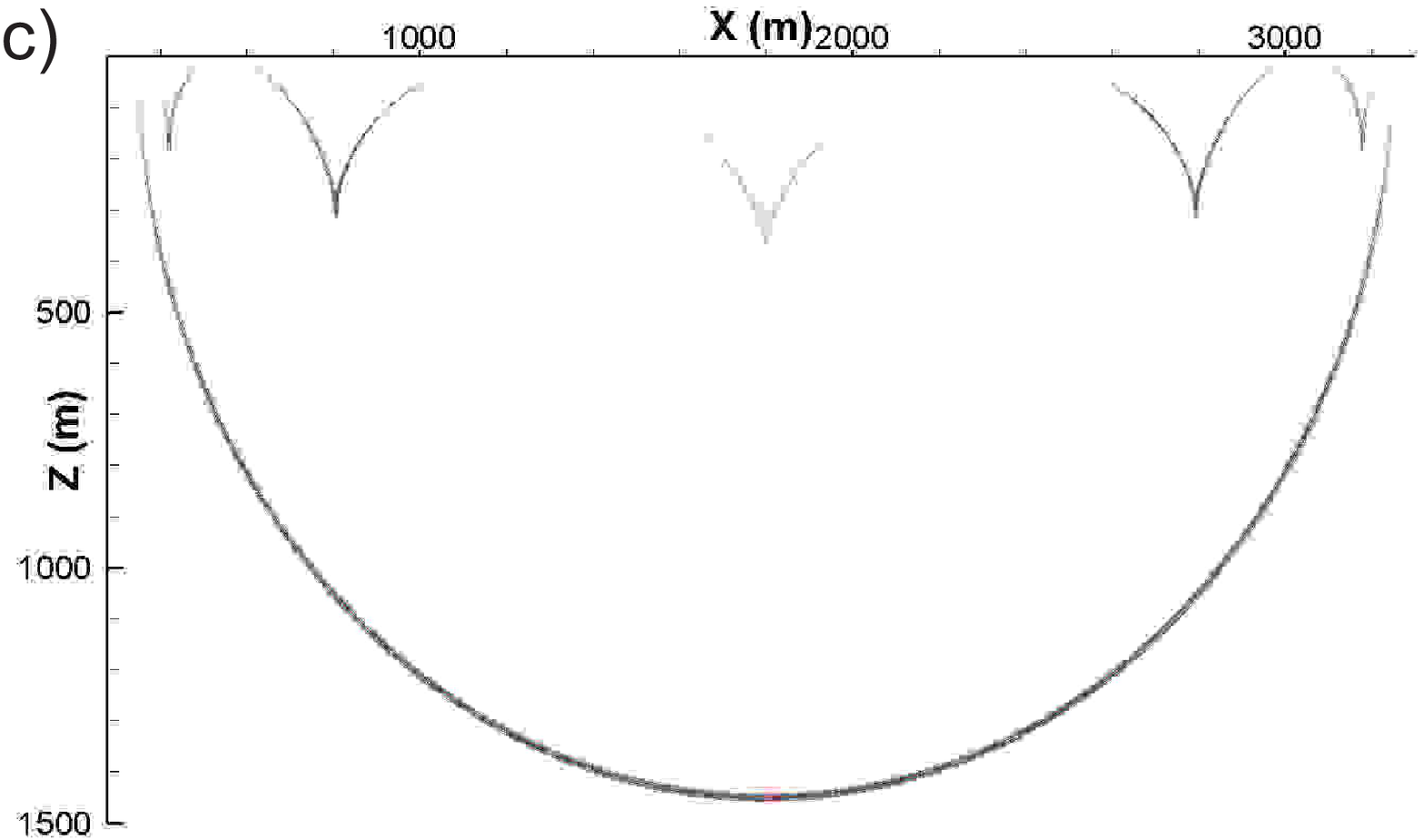}
        \end{subfigure}%
        \caption{(a) The 2D Sigsbee2A salt model and (b) zero-offset section.}
                \label{fig:sigsbee_model}
\end{figure}

For Sigsbee2A model \cite{Paffenholz} (Fig.~\ref{fig:sigsbee_model}a), theoretical seismograms (Fig.~\ref{fig:sigsbee_model}b) were calculated using the algorithm of explosive boundaries implemented in the Madagascar package \cite{Fomel2013}. For setting the boundary condition, the function for zero-offset section$\left.u(x,z,t)\right|_{z=0}=g(x,t)$ was expanded in series (\ref{series_lag}) for $n=3500$, the parameter  $\eta=300$ for $t\in[0,12]\;s$. The calculations were done on the meshes with the steps $h_{x,z}=12.5\;m$ for the Richardson method  and $h_{x}=12.5\;m$, $h_{z}=4.16\;m$ for the methods AM5-I5, PC5-I5.

The velocity model Sigsbee is not smooth, therefore preliminary smoothing (\ref{smooth}) for the stability of the algorithms  AM5-I5 and PC5-I5 is required. If for the model Syncline a single smoothing to provide the stability is needed, for the model Sigsbee the smoothing procedure was applied three times.  For the Richardson method the smoothing procedure of the velocity model because of a higher inner dissipation and stability, was not carried out.

\begin{figure}[!h]
\centering
                \includegraphics[width=\textwidth]{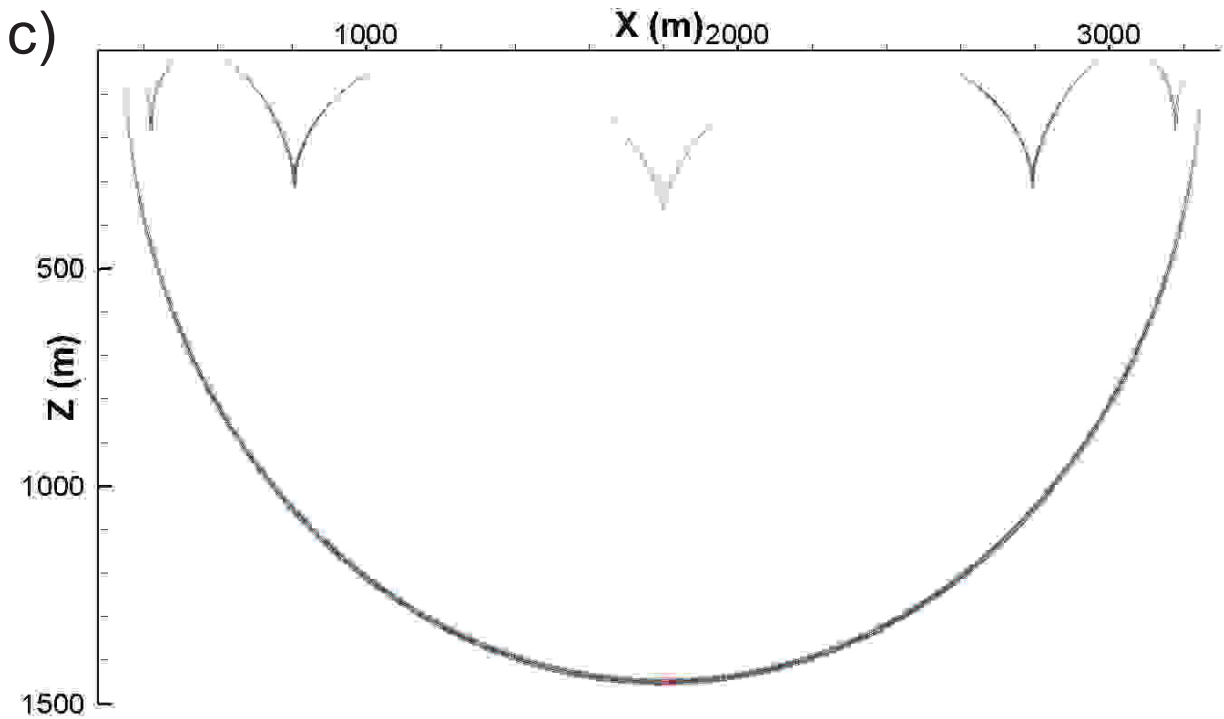}

\caption{Snapshots for the wave field at $t=12\;s$ for the model and the zero-offset section in Fig.~\ref{fig:sigsbee_model} for the Richardson method with $h_{x,z}=12.5$ m.}
\label{fig:sigsbee_rich}
 \end{figure}

The simulation results for the Richardson method are presented in Fig.~\ref{fig:sigsbee_rich}; the results obtained for the methods AM5-I5, PC5-I5 are given in Figure 9, because for these methods there is no difference in solution. The Richardson method with the total number of nodes of the meshes  $\Omega_1,\Omega_2$  equal to the number of nodes of the mesh $\Omega$ for the methods AM5-I5, PC5-I5 possesses less accuracy for amplitudes, which is due to a stronger numerical dissipation. All the three algorithms have demonstrated the stability, however the multistep procedures of the predictor-corrector type are more efficient.

 \begin{figure}[!h]
 \centering
                \includegraphics[width=\textwidth]{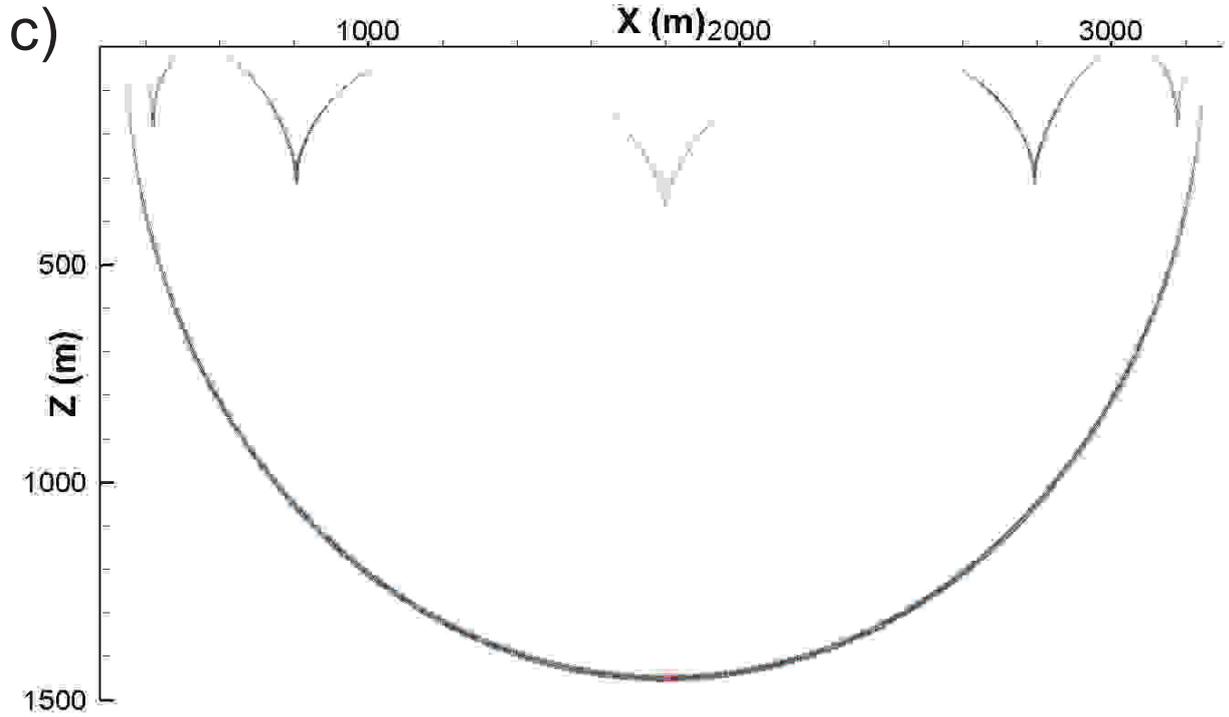}
\caption{Snapshots for the wave field at $t=12\;s$ for the model and the zero-offset section (Fig.~\ref{fig:sigsbee_model}) for the AM5-I5 and PC5-I5 methods with $h_x=12.5$ m and $h_z=4.16$ m.}
\label{fig:sigsbee_adams}
\end{figure}
\section{Conclusion}
As a result of the conducted study it appeared possible to determine the reasons of the numerical instability for the schemes of a high-order of accuracy in solving the OWWE and to propose the stabilizing procedure based on the spline-filtration. This allowed us to implement the Adams multistep schemes and the predictor-corrector method accurate to the fifth order. A combination of the Laguerre transform and the predictor-corrector method reduces the original problem to solving the SLAEs with symmetrical well-conditioned real matrices, which is one of the main advantages of the methods proposed as compared to the classical Fourier approach. In the 2D case, direct algorithms for solving SLAEs are applied, which are less sensitive to the above-mentioned properties of matrices. But for the 3D problems, conditioning of matrices is defined by the rate of convergence of iterative procedures and hence the total calculation time.

In spite of the fact that the Richardson extrapolation procedure demands is more computationally expensive costs than the Adams multistep methods, one should not completely reject its application. First, the Richardson method possesses a greater numerical dissipation, which in many cases makes possible to calculate inhomogeneous velocity models without preliminary smoothing. Also, supplementary stability will not be redundant when considering the OWWE for an elastic model. Second, the Richardson method can be used to calculate initial values for multistep methods that are not self-starting. As a rule, for solving this problem the Runge-Kutta type schemes are employed, but within the Laguerre method such schemes do not provide the required accuracy due to a strong numerical dissipation.

A combination of the spline-filtration, the Adams multistep methods and the Laguerre transform is mutually complementary. Experimentally it was verified that the change of the Adams methods for the  backward difference schemes does not provide the stability of calculation with the help of proposed stabilizing procedures, whereas the change of the Laguerre transform for the Fourier transform with respect to time makes the spline-filtration unreasonable. This is because in this case the solution for each harmonic is independently determined by initial conditions on the daily surface. On the contrary, the matter of coefficients of the Laguerre series and their recurrent dependence make possible to delicately remove unstable components of the wave field not additing new numerical artifacts. Thus, the considered ways of decreasing computer costs make the proposed methods of solving the OWWE to be promising for the calculation of applied problems.
\section{Acknowledgments}
The calculation was performed on the supercomputers of the Moscow and the Novosibirsk State Universities and of the Siberian supercomputer center.
The work was supported by a grant from the Russian Ministry of Education no. MK-152.2017.5.
\newpage
\bibliography{base}
\end{document}